\documentclass{article}
\usepackage{amsmath,mathrsfs,amssymb,amsthm,amscd}
\usepackage[]{fontenc}
\usepackage[all]{xy}
\usepackage{hyperref}
\usepackage{graphics}
\usepackage{a4wide}

\newcounter{EQNR}
\renewcommand{\contentsname}{Contents \\
{\footnotesize\normalfont(A table of contents should normally not be included)}}

\begin{document}

\title{Effective bounds for Faltings's delta function}
\author{\emph{Dedicated to Christophe Soul\'e at his sixtieth birthday}}
\date{Jay Jorgenson\footnote{The first author acknowledges support from numerous NSF
and PSC-CUNY grants.}\mbox{ } and J\"urg Kramer\footnote{The second author acknowledges
support from the DFG Graduate School \emph{Berlin Mathematical School} and from the
DFG International Research Training Group \emph{Moduli and Automorphic Forms}.}}

\maketitle

\section{Introduction}

In \cite{Arakelov1} and \cite{Arakelov2}, S.~J.~Arakelov introduced Green's functions on
compact Riemann surfaces in order to define an intersection theory on arithmetic surfaces,
thus initiating a far-reaching mathematical program which bears his name. G.~Faltings
extended the pioneering work of Arakelov in \cite{Faltings} by defining metrics on
determinant line bundles arising from the cohomology of algebraic curves, from which
he derived arithmetic versions of the Riemann-Roch theorem, Noether's formula, and
the Hodge index theorem. Although Faltings employed the classical Riemann theta
function to define metrics on these determinant line bundles, he does refer to the
emerging idea of D.~Quillen to use Ray--Singer analytic torsion to define the metrics
on these determinant line bundles as being ``more direct''.  One of the many aspects of
the mathematical legacy of Christophe Soul\'e is the central role he played in developing
higher dimensional Arakelov theory where Quillen metrics are fully utilized; see \cite
{Soule2} and the references therein.

From Faltings's theory \cite{Faltings}, there appears a naturally defined analytic quantity
associated to any compact Riemann surface $X$. The new invariant in \cite{Faltings} became 
known as Faltings's delta function, which we denote by $\delta_{\mathrm{Fal}}(X)$. Many
of the fundamental arithmetic theorems and formulas in \cite{Faltings}, such as those listed 
above, amount to statements which involve $\delta_{\mathrm{Fal}}(X)$. By comparing the
Riemann-Roch theorem from \cite{Faltings} and the arithmetic Riemann-Roch theorem,
Soul\'e expressed in \cite{Soule1} the Faltings's delta function in terms of the analytic torsion
of the trivial line bundle on $X$ when given the Arakelov metric; see equation \eqref{2}
below, as well as \cite{Soule1} and, more recently, \cite{Wentworth2}.  

The Polyakov formula allows one to relate values of the analytic torsion for conformally
equivalent metrics. As a result, one can use Soul\'e's formula for the Faltings's delta function
and obtain an identity which expresses $\delta_{\mathrm{Fal}}(X)$ in terms of the hyperbolic
geometry of $X$; see Theorem \ref{3.2} below, which comes from \cite{JK4}. In \cite{JK4},
we used the relation between the Faltings's delta function and the hyperbolic geometry in
order to study $\delta_{\mathrm{Fal}}(X)$ through covers. As an arithmetic application of
the analytic bounds obtained for $\delta_{\mathrm{Fal}}(X)$, we derived in \cite{JK4} an
improved estimate for the Faltings height of the Jacobian of the modular curve $X_{0}(N)$
for square-free $N$ which is not divisible by $6$.

After the completion of \cite{JK4}, A.~N.~Parshin posed the following question to the second
named author: Can one derive an \emph{effective} bound for the Faltings's delta function
$\delta_{\mathrm{Fal}}(X)$ in terms of basic information associated to the hyperbolic
geometry of $X$? The purpose of the present article is to provide an affirmative answer
to Parshin's question. More specifically, our main results, given in Theorem~\ref{6.1} and 
Corollaries~\ref{6.3}, \ref{6.4}, explicitly bound $\delta_{\mathrm{Fal}}(X)$ in the case
when $X$ is a finite degree covering of a compact Riemann surface $X_{0}$ of genus
bigger than $1$, where the bound for $\delta_{\mathrm{Fal}}(X)$ is effectively computable
once knowing the genera of $X_{0}$ and $X$, the smallest non-zero eigenvalues of the
hyperbolic Laplacian acting on $X_{0}$ and $X$, and the length of the shortest closed
geodesic on $X_{0}$ (as well as some ramification data in case the covering is ramified).

An important ingredient in the analysis of the present paper is the algorithm from \cite{FJK},
which provides effective means by which one can bound the Huber constant on $X$, a
quantity associated to the error term in the prime geodesic theorem; see \cite{FJK} and the
references therein. As with the main result in \cite{FJK}, it is possible that the effective bound
we obtain here may not be optimal, perhaps even far from it. However, the existence of an
effective bound for the Faltings's delta function $\delta_{\mathrm{Fal}}(X)$, albeit a sub-optimal
bound, may be a tool by which one can further investigate the application of Arakelov theory
to diophantine problems, as originally intended.

The paper is organized as follows: After recalling basic notations in section~2, we express
Faltings's delta function $\delta_{\mathrm{Fal}}(X)$ in hyperbolic terms of $X$ in section~3.
Section~4 is devoted to derive effective bounds for the ratio $\mu_{\mathrm{can}}(z)/\mu_
{\mathrm{hyp}}(z)$ of the canonical by the hyperbolic metric on $X$ and section~5 gives
effective bounds for the Huber constant $C_{\mathrm{Hub},X}$ on $X$. In section~6, we
combine the results of the sections~3, 4, 5 to derive effective bounds for $\delta_{\mathrm
{Fal}}(X)$. In section~7, we discuss an application of our results to an idea of A.~N.~Parshin
for an attempt giving effective bounds for the height of rational points on smooth projective
curves defined over number fields.

\emph{Acknowledgements:} We would like to use this opportunity to thank Christophe
Soul\'e for having introduced us into the theory of arithmetic intersections by generously
sharing his broad knowledge and deep insights on the subject with us. Furthermore, we
would like to thank Alexei Parshin for his interest in our results and for having pointed
out to us an application to his work. Finally, we would like to thank the referee for some
of his/her comments.

\section{Basic notations}

\begin{nn}\label{2.1}
\textbf{Hyperbolic and canonical metrics.} In this note $X$ will denote a compact Riemann
surface of genus $g_{X}>1$. By the uniformization theorem, $X$ is isomorphic to the quotient
space $\Gamma\backslash\mathbb{H}$, where $\Gamma$ is a cocompact and torsionfree
Fuchsian subgroup of the first kind of $\mathrm{PSL}_{2}(\mathbb{R})$ acting by fractional
linear transformations on the upper half-plane $\mathbb{H}=\{z\in\mathbb{C}\,\vert\,z=x+iy,\,
y>0\}$. In the sequel, we will identify $X$ locally with its universal cover $\mathbb{H}$. 

We denote by $\mu_{\mathrm{hyp}}$ the $(1,1)$-form corresponding to the hyperbolic metric on
$X$, which is compatible with the complex structure of $X$ and has constant negative curvature
equal to $-1$. Locally, we have
\begin{align*}
\mu_{\mathrm{hyp}}(z)=\frac{i}{2}\cdot\frac{\mathrm{d}z\wedge\mathrm{d}\overline{z}}{\mathrm
{Im}(z)^{2}}=\frac{\mathrm{d}x\wedge\mathrm{d}y}{y^{2}}\,.
\end{align*}
We write $\mathrm{vol}_{\mathrm{hyp}}(X)$ for the hyperbolic volume of $X$; recall that $\mathrm
{vol}_{\mathrm{hyp}}(X)$ is given by $4\pi(g_{X}-1)$. By $\mu_{\mathrm {shyp}}$, we denote the
$(1,1)$-form corresponding to the rescaled hyperbolic metric, which measures the volume of $X$
to be $1$. We write $\mathrm{dist}_{\mathrm{hyp}}(z,w)$ for the hyperbolic distance between two
points $z,w\in\mathbb{H}$. We recall the formula
\begin{align*}
\mathrm{dist}_{\mathrm{hyp}}(z,w)=\cosh^{-1}\bigg(1+\frac{\vert z-w\vert^{2}}{2\,\mathrm{Im}(z)
\mathrm{Im}(w)}\bigg).
\end{align*}
We denote the hyperbolic Laplacian on $X$ by $\Delta_{\mathrm{hyp}}$; locally, we have
\begin{align*}
\Delta_{\mathrm{hyp}}=-y^{2}\bigg(\frac{\partial^{2}}{\partial x^{2}}+\frac{\partial^{2}}{\partial y^
{2}}\bigg).
\end{align*}
The discrete spectrum of $\Delta_{\mathrm{hyp}}$ is given by the increasing sequence of
eigenvalues 
\begin{align*}
0=\lambda_{X,0}<\lambda_{X,1}\leq\lambda_{X,2}\leq\dots
\end{align*}
The $(1,1)$-form $\mu_{\mathrm{can}}$ associated to the canonical metric is defined as follows.
Let $\{\omega_{1},\ldots,\omega_{g_{X}}\}$ denote an orthonormal basis of the space $\Gamma
(X,\Omega_{X}^{1})$ of holomorphic $1$-forms on $X$. Then, $\mu_{\mathrm{can}}$ is locally 
given by
\begin{align*}
\mu_{\mathrm{can}}(z)=\frac{1}{g_{X}}\cdot\frac{i}{2}\,\sum_{j=1}^{g_{X}}\omega_{j}(z)\wedge
\overline{\omega}_{j}(z).
\end{align*}
We recall that the Arakelov metric on $X$ is induced by means of the residual canonical metric
$\Vert\cdot\Vert_{\mathrm{Ar}}$ on $\Omega_{X}^{1}$, which turns the residue map into an
isometry.
\end{nn}

\begin{nn}\label{2.2}
\textbf{Hyperbolic heat kernel for functions.} The hyperbolic heat kernel $K_{\mathbb{H}}(t;z,w)$
($t\in\mathbb{R}_{>0}$; $z,w\in\mathbb{H}$) for functions on $\mathbb{H}$ is given by the formula
\begin{align*}
K_{\mathbb{H}}(t;z,w):=K_{\mathbb{H}}(t;\rho):=\frac{\sqrt{2}e^{-t/4}}{(4\pi t)^{3/2}}\int\limits_{\rho}^
{\infty}\frac{re^{-r^{2}/(4t)}}{\sqrt{\cosh(r)-\cosh(\rho)}}\,\mathrm{d}r\,,
\end{align*}
where $\rho=\mathrm{dist}_{\mathrm{hyp}}(z,w)$. The hyperbolic heat kernel $K_{X}(t;z,w)$ ($t\in
\mathbb{R}_{>0}$; $z,w\in X$) for functions on $X$ is obtained by averaging over the elements of
$\Gamma$, namely
\begin{align*}
K_{X}(t;z,w):=\sum_{\gamma\in\Gamma}K_{\mathbb{H}}(t;z,\gamma w).
\end{align*}
The heat kernel $K_{X}(t;z,w)$ satisfies the equations
\begin{align*}
&\bigg(\frac{\partial}{\partial t}+\Delta_{\mathrm{hyp},z}\bigg)K_{X}(t;z,w)=0\qquad\qquad\,(z,w
\in X), \\
&\lim_{t\rightarrow 0}\int\limits_{X}K_{X}(t;z,w)\,f(w)\,\mu_{\mathrm{hyp}}(w)=f(z)\quad(z\in X)
\end{align*}
for all $C^{\infty}$-functions $f$ on $X$. As a shorthand, we use in the sequel the notation
\begin{align*}
HK_{X}(t;z):=\sum_{\substack{\gamma\in\Gamma\\\gamma\neq\mathrm{id}}}K_{\mathbb{H}}
(t;z,\gamma z).
\end{align*}
\end{nn}

\begin{nn}\label{2.3}
\textbf{Selberg zeta function.} Let $\mathcal{H}(\Gamma)$ denote the set of conjugacy classes
of primitive, hyperbolic elements in $\Gamma$. We denote by $\ell_{\gamma}$ the  hyperbolic
length of the closed geodesic determined by $\gamma\in\mathcal{H}(\Gamma)$ on $X$; it is
well-known that the equality
\begin{align*}
\vert\mathrm{tr}(\gamma)\vert=2\cosh(\ell_{\gamma}/2)
\end{align*}
holds.

For $s\in\mathbb{C}$, $\mathrm{Re}(s)>1$, the Selberg zeta function $Z_{X}(s)$ associated
to $X$ is defined via the Euler product expansion
\begin{align*}
Z_{X}(s):=\prod_{\gamma\in\mathcal{H}(\Gamma)}Z_{\gamma}(s),
\end{align*}
where the local factors $Z_{\gamma}(s)$ are given by
\begin{align*}
Z_{\gamma}(s):=\prod_{n=0}^{\infty}\big(1-e^{-(s+n)\ell_{\gamma}}\big).
\end{align*}
The Selberg zeta function $Z_{X}(s)$ is known to have a meromorphic continuation to all of
$\mathbb{C}$ with zeros and poles characterized by the spectral theory of the hyperbolic
Laplacian; furthermore, $Z_{X}(s)$ satisfies a functional equation. For our purposes, it suffices
to know that the Selberg zeta function $Z_{X}(s)$ has a simple zero at $s=1$, so that the
quantity
\begin{align*}
\lim_{s\rightarrow 1}\bigg(\frac{Z_{X}'}{Z_{X}}(s)-\frac{1}{s-1}\bigg)
\end{align*}
is well-defined.
\end{nn}

\begin{nn}\label{2.4}
\textbf{Prime geodesic theorem.} For any small eigenvalue $\lambda_{X,j}\in[0,1/4)$, we define
\begin{align*}
s_{X,j}:=\frac{1}{2}+\sqrt{\frac{1}{4}-\lambda_{X,j}},
\end{align*}
and note that $1/2<s_{X,j}\leq 1$. For $u\in\mathbb{R}_{>1}$, we recall the prime geodesic
counting function
\begin{align*}
\pi_{X}(u):=\#\big\{\gamma\in\mathcal{H}(\Gamma)\,\vert\,e^{\ell_{\gamma}}<u\big\}.
\end{align*}
Introducing the logarithmic  integral
\begin{align*}
\mathrm{li}(u):=\int\limits_{2}^{u}\frac{\mathrm{d}\xi}{\log(\xi)}\,,
\end{align*}
the prime geodesic theorem states
\begin{align}
\label{1}
\pi_{X}(u)=\sum_{0\leq\lambda_{X,j}<1/4}\mathrm{li}(u^{s_{X,j}})+O_{X}\big(u^{3/4}\log(u)^
{-1/2}\big)
\end{align}
for $u>1$, where the implied constant for all $u>1$, not just asymptotically, depends solely on
$X$. We call the infimum of all possible implied constants the Huber constant and denote it by
$C_{\mathrm{Hub},X}$.
\end{nn}

\section{Faltings's delta function in hyperbolic terms}

\begin{nn}\label{3.1} 
\textbf{Faltings's delta function.} Faltings's delta function $\delta_{\mathrm{Fal}}(X)$ was introduced in
\cite{Faltings}, where also some of its basic properties were given. In \cite{Jorgenson}, Faltings's delta
function is expressed in terms of Riemann theta functions, and its asymptotic behavior is investigated;
see also \cite{Wentworth1}. As a by-product of the analytic part of the arithmetic Riemann-Roch theorem
for arithmetic surfaces, C.~Soul\'e has shown in \cite{Soule1} that
\begin{align}
\label{2}
\delta_{\mathrm{Fal}}(X)=-6D_{\mathrm{Ar}}(X)+a(g_{X})\,,
\end{align}
where
\begin{align*}
D_{\mathrm{Ar}}(X):=\log\bigg(\frac{\det^{*}(\Delta_{\mathrm{Ar}})}{\mathrm{vol}_{\mathrm{Ar}}
(X)}\bigg)
\end{align*}
with $\det^{*}(\Delta_{\mathrm{Ar}})$ the regularized determinant  of the Laplacian, $\mathrm{vol}_
{\mathrm{Ar}}(X)$ the volume with respect to the Arakelov metric $\Vert\cdot\Vert_{\mathrm{Ar}}$,
and
\begin{align*}
a(g_{X}):=-2g_{X}\log(\pi)+4g_{X}\log(2)+(g_{X}-1)(-24\zeta_{\mathbb{Q}}'(-1)+1).
\end{align*}
\end{nn}

It has been shown in \cite{JK4} how Faltings's delta function can be expressed solely in hyperbolic
terms. Theorem~3.8 therein states:

\begin{nn}\label{3.2}
\textbf{Theorem.} \emph{For $X$ with genus $g_{X}>1$, let
\begin{align*}
F(z):=\int\limits_{0}^{\infty}\bigg(HK_{X}(t;z)-\frac{1}{\mathrm{vol}_{\mathrm{hyp}}(X)}\bigg)\mathrm{d}t.
\end{align*}
Then, we have
\begin{align}
\notag
&\delta_{\mathrm{Fal}}(X)= \\
\label{3}
&2\pi\bigg(1-\frac{1}{g_{X}}\bigg)\int\limits_{X}F(z)\Delta_{\mathrm{hyp}}F(z)\mu_{\mathrm{hyp}}(z)-
6\log\big(Z_{X}'(1)\big)+2\lim_{s\rightarrow 1}\bigg(\frac{Z_{X}'}{Z_{X}}(s)-\frac{1}{s-1}\bigg)+c(g_{X}),
\end{align}
where
\begin{align*}
c(g_{X}):=&\,\,a(g_{X})-6b(g_{X})+2(g_{X}-1)\log(4)+6\log\big(\mathrm{vol}_{\mathrm{hyp}}(X)\big)-2 \\
=&\,\,2g_{X}\big(-24\zeta_{\mathbb{Q}}'(-1)-4\log(\pi)+\log(2)+2\big)+6\log\big(\mathrm{vol}_{\mathrm
{hyp}}(X)\big)+ \\
&\,\,\big(48\zeta_{\mathbb{Q}}'(-1)+6\log(2\pi)-2\log(4)-6\big)
\end{align*}
with $a(g_{X})$ as above and $b(g_{X})$ given by}
\begin{align*}
b(g_{X}):=(g_{X}-1)\big(4\zeta_{\mathbb{Q}}'(-1)-1/2+\log(2\pi)\big).
\end{align*}
\begin{proof}
The proof is given in \cite{JK4}. Here we present only a short outline of the proof, which consists
of the following three main ingredients: \\
(i) One starts by using the Polyakov formula to relate the regularized determinants with respect to
the Arakelov and the hyperbolic metric, namely
\begin{align*}
D_{\mathrm{Ar}}(X)=D_{\mathrm{hyp}}(X)+\frac{g_{X}-1}{6}\int\limits_{X}\phi_{\mathrm{Ar}}(z)\big
(\mu_{\mathrm{can}}(z)+\mu_{\mathrm{hyp}}(z)\big),
\end{align*}
where $\phi_{\mathrm{Ar}}(z)$ is the conformal factor describing the change from the Arakelov to
the hyperbolic metric. \\
(ii) In a second step, one uses the result \cite{Sarnak} by P.~Sarnak describing the hyperbolic
regularized determinant in terms of the Selberg zeta function, namely
\begin{align*}
D_{\mathrm{hyp}}(X)=\log\bigg(\frac{Z_{X}'(1)}{\mathrm{vol}_{\mathrm{hyp}}(X)}\bigg)+b(g_{X}).
\end{align*}
(iii) In order to express the conformal factor $\phi_{\mathrm{Ar}}(z)$ and the canonical metric
form $\mu_{\mathrm{can}}(z)$ in hyperbolic terms, we make use of the fundamental relation
\begin{align}
\label{4}
\mu_{\mathrm{can}}(z)=\mu_{\mathrm{shyp}}(z)+\frac{1}{2g_{X}}\Bigg(\int\limits_{0}^{\infty}\Delta_
{\mathrm{hyp}}K_{X}(t;z)\,\mathrm{d}t\Bigg)\mu_{\mathrm{hyp}}(z),
\end{align}
which has been proven in Appendix~1 of \cite{JK4}.
\end{proof}
\end{nn}

\begin{nn}\label{3.3}
\textbf{Remark.} We note that formula \eqref{4} has meanwhile been generalized to cofinite
Fuchsian subgroups of the first kind of $\mathrm{PSL}_{2}(\mathbb{R})$ without torsion elements
in \cite{JK6}, and, as a relation of $(1,1)$-currents, to cofinite Fuchsian subgroups of the first kind
of $\mathrm{PSL}_{2}(\mathbb{R})$ allowing torsion elements in~\cite{Anil}. \\
\end{nn}

Based on formula \eqref{3}, the following bound can be derived for $\delta_{\mathrm{Fal}}(X)$ in
terms of basic hyperbolic invariants of $X$. For this we introduce the following notations
\begin{align*}
&\lambda_{X}:=\frac{1}{2}\min\bigg\{\lambda_{X,1},\frac{7}{64}\bigg\}, \\
&N_{\mathrm{ev},X}^{[0,1/4)}:=\#\big\{\lambda_{X,j}\,\vert\,0\leq\lambda_{X,j}<1/4\big\}, \\[1mm]
&N_{\mathrm{geo},X}^{(0,5)}:=\#\big\{\gamma\in\mathcal{H}(\Gamma)\,\vert\,0<\ell_{\gamma}<5
\big\}, \\
&S_{X}:=\sup\limits_{z\in X}\bigg(\frac{\mu_{\mathrm{can}}(z)}{\mu_{\mathrm{shyp}}(z)}\bigg),
\end{align*}
where $\lambda_{X,1}$ is the smallest non-zero eigenvalue of $\Delta_{\mathrm{hyp}}$, and we
recall that $\ell_{X}$ denotes the length of the shortest closed geodesic on $X$ and $C_{\mathrm
{Hub},X}$ is the Huber constant introduced in subsection~\ref{2.4}.

\begin{nn}\label{3.4}
\textbf{Corollary.} \emph{With the above notations, we have the bound
\begin{align*}
\delta_{\mathrm{Fal}}(X)\leq D_{1}\bigg(g_{X}+\frac{1}{\lambda_{X}}\big(g_{X}(S_{X}+1)^{2}+C_
{\mathrm{Hub},X}+N_{\mathrm{ev},X}^{[0,1/4)}\big)+\bigg(1+\frac{1}{\ell_{X}}\bigg)N_{\mathrm
{geo},X}^{(0,5)}\bigg)
\end{align*}
with an absolute constant $D_{1}>0$, which can be taken to be $10^{3}$.}
\begin{proof}
The proof is straightforward using Theorem~\ref{3.2} in combination with the estimates given in
Propositions~4.1, 4.2, 4.3, and Lemma~4.4 in \cite{JK4}. For the convenience of the reader, we
give now a more detailed derivation of the proof.

Using Proposition~4.1 of \cite{JK4} in combination with the inequalities $\lambda_{X,1}\geq
\lambda_{X}$ and $\mathrm{vol}_{\mathrm{hyp}}(X)\leq 4\pi g_{X}$, the integral in \eqref{3}
can be bounded as
\begin{align}
\label{5}
\bigg\vert\int\limits_{X}F(z)\Delta_{\mathrm{hyp}}F(z)\mu_{\mathrm{hyp}}(z)\bigg\vert\leq\frac{(S_
{X}+1)^{2}\,\mathrm{vol}_{\mathrm{hyp}}(X)}{\lambda_{X,1}}\leq\frac{4\pi g_{X}}{\lambda_{X}}(S_
{X}+1)^{2}.
\end{align}
In order to bound the absolute value of the second summand in \eqref{3}, we first observe that
we have to take the second bound in Proposition~4.3 of \cite{JK4}, since the first one being
logarithmic in $g_{X}$ is too small; choosing $\varepsilon=\lambda_{X}$, we obtain
\begin{align*}
\big\vert \log\big(Z_{X}'(1)\big)\big\vert\leq -\sum_{\substack{\gamma\in\mathcal{H}(\Gamma)\\
\ell_{\gamma}<5}}\log\big(Z_{\gamma}(1)\big)+12\bigg(5+\frac{1}{\lambda_{X}}\bigg)\big(C_
{\mathrm{Hub},X}+N_{\mathrm{ev},X}^{[0,1/4)}+1\big).
\end{align*}
Using Lemma~4.4 (i) of \cite{JK4}, we derive from this the bound
\begin{align}
\notag
\big\vert \log\big(Z_{X}'(1)\big)\big\vert&\leq\sum_{\substack{\gamma\in\mathcal{H}(\Gamma)\\
\ell_{\gamma}<5}}\frac{\pi^{2}}{6\ell_{\gamma}}+\frac{72}{\lambda_{X}}\big(C_{\mathrm{Hub},
X}+N_{\mathrm{ev},X}^{[0,1/4)}+1\big) \\
\label{6}
&\leq\,\frac{\pi^{2}}{6\ell_{X}}N_{\mathrm{geo},X}^{(0,5)}+\frac{144}{\lambda_{X}}\big(C_{\mathrm
{Hub},X}+N_{\mathrm{ev},X}^{[0,1/4)}\big).
\end{align}
Finally, in order to bound the absolute value of the third summand in \eqref{3}, we again observe
that we have to take the second bound in Proposition~4.2 of \cite{JK4}, since the first one being
logarithmic in $g_{X}$ is too small; choosing again $\varepsilon=\lambda_{X}$, we obtain
\begin{align*}
\bigg\vert\lim_{s\rightarrow 1}\bigg(\frac{Z_{X}'}{Z_{X}}(s)-\frac{1}{s-1}\bigg)\bigg\vert\leq\sum_
{\substack{\gamma\in\mathcal{H}(\Gamma)\\\ell_{\gamma}<5}}\frac{Z_{\gamma}'}{Z_{\gamma}}
(1)+\frac{6}{\lambda_{X}}\big(C_{\mathrm{Hub},X}+N_{\mathrm{ev},X}^{[0,1/4)}\big)+2\,.
\end{align*}
Using Lemma~4.4 (ii) of \cite{JK4}, we derive from this the bound
\begin{align}
\notag
\bigg\vert\lim_{s\rightarrow 1}\bigg(\frac{Z_{X}'}{Z_{X}}(s)-\frac{1}{s-1}\bigg)\bigg\vert&\leq\sum_
{\substack{\gamma\in\mathcal{H}(\Gamma)\\\ell_{\gamma}<5}}\bigg(3+\log\bigg(\frac{1}{\ell_
{\gamma}}\bigg)\bigg)+\frac{6}{\lambda_{X}}\big(C_{\mathrm{Hub},X}+N_{\mathrm{ev},X}^{[0,
1/4)}\big)+2 \\
\notag
&\leq\sum_{\substack{\gamma\in\mathcal{H}(\Gamma)\\\ell_{\gamma}<5}}\bigg(3+\frac{1}{\ell_
{\gamma}}\bigg)+\frac{6}{\lambda_{X}}\big(C_{\mathrm{Hub},X}+N_{\mathrm{ev},X}^{[0,1/4)}
\big)+2 \\
\label{7}
&\leq\bigg(3+\frac{1}{\ell_{X}}\bigg)N_{\mathrm{geo},X}^{(0,5)}+\frac{6}{\lambda_{X}}\big(C_
{\mathrm{Hub},X}+N_{\mathrm{ev},X}^{[0,1/4)}\big)+2\,.
\end{align}
The quantity $c(g_{X})$ in \eqref{3} is easily bounded as
\begin{align}
\label{8}
c(g_{X})\leq 11g_{X}+10.
\end{align}
Adding up the bounds \eqref{5}--\eqref{8}, using that $g_{X}>1$, and by crudely estimating the
arising integral constants by $D_{1}=10^{3}$, yields the claimed bound. Note that, estimating
more rigorously, $D_{1}$ can in fact be taken to be $876$.
\end{proof}
\end{nn}

\section{Effective bounds for the sup-norm}

\begin{nn}\label{4.1} 
\textbf{Hyperbolic heat kernel for forms.} In addition to the hyperbolic heat kernel on $\mathbb{H}$, resp.
$X$, introduced in subsection \ref{2.2}, we also need the hyperbolic heat kernel for forms of weight $1$
on $\mathbb{H}$, resp. $X$. The hyperbolic heat kernel for forms of weight $1$ on $\mathbb{H}$ is
defined as in \cite{JK3}, namely we have
\begin{align*}
K_{\mathbb{H}}^{(1)}(t;z,w):=K_{\mathbb{H}}^{(1)}(t;\rho):=\frac{\sqrt{2}e^{-t/4}}{(4\pi t)^{3/2}}\int\limits_
{\rho}^{\infty}\frac{re^{-r^{2}/(4t)}}{\sqrt{\cosh(r)-\cosh(\rho)}}\,T_{2}\bigg(\frac{\cosh(r/2)}{\cosh(\rho/2)}
\bigg)\,\mathrm{d}r\,,
\end{align*}
where $T_{2}$ is the Chebyshev polynomial given by $T_{2}(r):=2r^{2}-1$. The hyperbolic heat kernel
for forms of weight $1$ on $X$ on the diagonal is then given as
\begin{align*}
K_{X}^{(1)}(t;z):=\sum_{\gamma\in\Gamma}c(\gamma;z)\,K_{\mathbb{H}}^{(1)}(t;z,\gamma z),
\end{align*}
where $c(\gamma,z)$ for $\gamma=\big(\begin{smallmatrix}a&b\\c&d\end{smallmatrix}\big)$ is defined
as
\begin{align*}
c(\gamma,z):=\frac{c\bar{z}+d}{cz+d}\cdot\frac{z-\gamma\bar{z}}{\gamma z-\bar{z}}\,.
\end{align*}
We note that
$\vert c(\gamma,z)\vert=1$. From \cite{JK3}, we recall the crucial relation
\begin{align}
\label{9}
\lim\limits_{t\rightarrow\infty}K_{X}^{(1)}(t;z)=\frac{g_{X}\mu_{\mathrm{can}}(z)}{\mu_{\mathrm{hyp}}
(z)}\,.
\end{align}
\end{nn}

\begin{nn}\label{4.2} 
\textbf{Lemma.} \emph{With the above notations, we have the bound
\begin{align}
\notag
K_{\mathbb{H}}^{(1)}(t;\rho)&\leq\frac{17\sqrt{2}\,e^{-t/4}}{(4\pi t)^{3/2}}\frac{(\rho+\log(4))e^{-\rho^{2}/
(4t)}}
{\sinh^{1/2}(\rho)}+\frac{4\,\sqrt{2}\,e^{-(\rho/(2\sqrt{t})+\sqrt{t}/2)^{2}}}{\pi^{3/2}\,\sqrt{t}} \\
\label{10}
&+\frac{4\,\sqrt{2}\,e^{-\rho}}{\pi^{3/2}}\int\limits_{\rho/(2\sqrt{t})-\sqrt{t}/2}^{\infty}e^{-r'^{2}}\,\mathrm{d}r'
\end{align}
for any $t>0$ and $\rho>0$.}
\begin{proof}
Starting with the defining formula
\begin{align*}
K_{\mathbb{H}}^{(1)}(t;\rho):=\frac{\sqrt{2}\,e^{-t/4}}{(4\pi t)^{3/2}}\int\limits_{\rho}^{\infty}\frac{re^{-r^{2}/
(4t)}}{\sqrt{\cosh(r)-\cosh(\rho)}}\,T_{2}\bigg(\frac{\cosh(r/2)}{\cosh(\rho/2)}\bigg)\,\mathrm{d}r\,,
\end{align*}
we decompose the integral under consideration as
\begin{align}
\label{11}
\int\limits_{\rho}^{\infty}\,\,\ldots\quad=\int\limits_{\rho}^{\rho+\log(4)}\ldots\quad+\int\limits_{\rho+\log
(4)}^{\infty}\ldots
\end{align}
We start by estimating the first integral on the right-hand side of \eqref{11}. Using the mean value
theorem for the function $\cosh(r)$ with $r\in[\rho,\rho+\log(4)]$, we obtain the bound
\begin{align*}
\cosh(r)-\cosh(\rho)=(r-\rho)\sinh(r_{*})\geq(r-\rho)\sinh(\rho),
\end{align*}
where $r_{*}\in[\rho,\rho+\log(4)]$. With this in mind, we have the estimate
\begin{align*}
&\int\limits_{\rho}^{\rho+\log(4)}\frac{re^{-r^{2}/(4t)}}{\sqrt{\cosh(r)-\cosh(\rho)}}\,T_{2}\bigg(\frac{\cosh
(r/2)}{\cosh(\rho/2)}\bigg)\,\mathrm{d}r\leq \\
&\frac{(\rho+\log(4))e^{-\rho^{2}/(4t)}}{\sinh^{1/2}(\rho)}T_{2}\bigg(\frac{\cosh((\rho+\log(4))/2)}{\cosh
(\rho/2)}\bigg)\int\limits_{\rho}^{\rho+\log(4)}(r-\rho)^{-1/2}\,\mathrm{d}r\leq \\
&\frac{2\log(4)^{1/2}(\rho+\log(4))e^{-\rho^{2}/(4t)}}{\sinh^{1/2}(\rho)}T_{2}\bigg(\frac{\cosh((\rho+\log
(4))/2)}{\cosh(\rho/2)}\bigg).
\end{align*}
Since, for any $r_{1},r_{2}\in\mathbb{R}_{>0}$, we have
\begin{align*}
\frac{\cosh(r_{1}+r_{2})}{\cosh(r_{1})}=\frac{\cosh(r_{1})\cosh(r_{2})}{\cosh(r_{1})}+\frac{\sinh(r_{1})
\sinh(r_{2})}{\cosh(r_{1})}\leq\cosh(r_{2})+\sinh(r_{2})=e^{r_{2}},
\end{align*}
we can estimate the Tshebyshev polynomial contribution as
\begin{align*}
T_{2}\bigg(\frac{\cosh((\rho+\log(4))/2)}{\cosh(\rho/2)}\bigg)\leq T_{2}\big(e^{\log(4)/2}\big)=7.
\end{align*}
In summary, we find the following bound for the integral in question
\begin{align}
\label{12}
&\int\limits_{\rho}^{\rho+\log(4)}\frac{re^{-r^{2}/(4t)}}{\sqrt{\cosh(r)-\cosh(\rho)}}\,T_{2}\bigg(\frac{\cosh
(r/2)}{\cosh(\rho/2)}\bigg)\,\mathrm{d}r\leq\frac{17(\rho+\log(4))e^{-\rho^{2}/(4t)}}{\sinh^{1/2}(\rho)}\,.
\end{align}
We now estimate the second integral on the right-hand side of \eqref{11}. Since $r\geq\rho+\log(4)$,
we have
\begin{align*}
\frac{\cosh(r)}{2}\geq\frac{\cosh(\rho+\log(4))}{2}\geq\frac{\cosh(\rho)\cosh(\log(4))}{2}\geq\cosh(\rho),
\end{align*}
whence
\begin{align*}
\cosh(r)-\cosh(\rho)\geq\frac{\cosh(r)}{2}\geq\frac{e^{r}}{4}\,.
\end{align*}
Therefore, using the estimate $T_{2}(r)\leq 2r^{2}$ in combination with 
\begin{align*}
\cosh(r/2)\leq e^{r/2}\quad\text{and}\quad\cosh(\rho/2)\geq\frac{e^{\rho/2}}{2}\,,
\end{align*}
we derive the bound
\begin{align}
\notag
&\int\limits_{\rho+\log(4)}^{\infty}\frac{re^{-r^{2}/(4t)}}{\sqrt{\cosh(r)-\cosh(\rho)}}\,T_{2}\bigg(\frac
{\cosh(r/2)}{\cosh(\rho/2)}\bigg)\,\mathrm{d}r\leq \\
\label{13}
&\int\limits_{\rho+\log(4)}^{\infty}\frac{2\,re^{-r^{2}/(4t)}}{e^{r/2}}\frac{8\,e^{r}}{e^{\rho}}\,\mathrm{d}
r=16\,e^{-\rho}\int\limits_{\rho+\log(4)}^{\infty}r\,e^{r/2}\,e^{-r^{2}/(4t)}\,\mathrm{d}r\,.
\end{align}
In order to complete the proof, we will further estimate the integral in \eqref{13}. Keeping in mind
that we finally have to multiply \eqref{13} by the factor $e^{-t/4}$, we estimate the quantity
\begin{align*}
&e^{-t/4}\int\limits_{\rho+\log(4)}^{\infty}r\,e^{r/2}\,e^{-r^{2}/(4t)}\,\mathrm{d}r\leq\int\limits_{\rho}^
{\infty}r\,e^{-(r/(2\sqrt{t})-\sqrt{t}/2)^{2}}\,\mathrm{d}r= \\
&2\,\sqrt{t}\int\limits_{\rho/(2\sqrt{t})-\sqrt{t}/2}^{\infty}\big(2\sqrt{t}r'+t\big)e^{-r'^{2}}\,\mathrm{d}r'=
2\,t\,e^{-(\rho/(2\sqrt{t})-\sqrt{t}/2)^{2}}+\,2\,t^{3/2}\int\limits_{\rho/(2\sqrt{t})-\sqrt{t}/2}^{\infty}e^{-r'^
{2}}\,\mathrm{d}r'\,.
\end{align*}
Multiplying by the remaining factor
\begin{align*}
\frac{16\,\sqrt{2}\,e^{-\rho}}{(4\pi t)^{3/2}}=\frac{2\,\sqrt{2}\,e^{-\rho}}{(\pi t)^{3/2}}\,,
\end{align*}
yields the following bound involving the second integral
\begin{align}
\notag
&\frac{\sqrt{2}\,e^{-t/4}}{(4\pi t)^{3/2}}\int\limits_{\rho+\log(4)}^{\infty}\frac{re^{-r^{2}/(4t)}}{\sqrt{\cosh
(r)-\cosh(\rho)}}\,T_{2}\bigg(\frac{\cosh(r/2)}{\cosh(\rho/2)}\bigg)\,\mathrm{d}r\leq \\
\label{14}
&\frac{4\,\sqrt{2}\,e^{-(\rho/(2\sqrt{t})+\sqrt{t}/2)^{2}}}{\pi^{3/2}\,\sqrt{t}}+\frac{4\,\sqrt{2}\,e^{-\rho}}
{\pi^{3/2}}\int\limits_{\rho/(2\sqrt{t})-\sqrt{t}/2}^{\infty}e^{-r'^{2}}\,\mathrm{d}r'\,.
\end{align}
Adding up the bounds \eqref{12} and \eqref{14} yields the claimed upper bound for $K_{\mathbb
{H}}^{(1)}(t;\rho)$.
\end{proof}
\end{nn}

\begin{nn}\label{4.3} 
\textbf{Lemma.} \emph{Let $X\longrightarrow X_{0}$ be an unramified covering of finite degree
with $X_{0}:=\Gamma_{0}\backslash\mathbb{H}$ a compact Riemann surface of genus $g_{0}>1$,
and let $\ell_{X_{0}}$ denote the length of the shortest closed geodesic on $X_{0}$. Then, the
quantity $S_{X}$ can be bounded as
\begin{align}
\label{15}
S_{X}\leq 4\pi\int\limits_{\ell_{X_{0}}/4}^{\infty}K_{\mathbb{H}}^{(1)}(t_{0};\rho)\,\frac{\sinh(\rho+\ell_
{X_{0}}/2)}{2\,\sinh^{2}(\ell_{X_{0}}/8)}\,\mathrm{d}\rho+4\pi\,K_{\mathbb{H}}^{(1)}(t_{0};\ell_{X_{0}}/4)
\bigg(\frac{\sinh^{2}(3\ell_{X_{0}}/8)}{\sinh^{2}(\ell_{X_{0}}/8)}-1\bigg)
\end{align}
for any $t_{0}>0$.}
\begin{proof}
From the spectral expansion, one immediately sees that the function $K_{X}^{(1)}(t;z)$ is monotone
decreasing in $t$. Using relation \eqref{9} together with the triangle inequality, we then obtain for
any $t_{0}>0$, the bound
\begin{align*}
\frac{g_{X}\mu_{\mathrm{can}}(z)}{\mu_{\mathrm{hyp}}(z)}\leq\sum_{\gamma\in\Gamma}K_
{\mathbb{H}}^{(1)}\big(t_{0};z,\gamma z\big)\leq\sum_{\gamma\in\Gamma_{0}}K_{\mathbb{H}}^
{(1)}\big(t_{0};z,\gamma z\big).
\end{align*}
Using the counting function
\begin{align*}
N_{X_{0}}(\rho;z):=\#\big\{\gamma\in\Gamma_{0}\,\vert\,\mathrm{dist}_{\mathrm{hyp}}(z,\gamma
z)<\rho\big\},
\end{align*}
we can express the latter bound in terms of the Stieltjes integral
\begin{align*}
\frac{g_{X}\mu_{\mathrm{can}}(z)}{\mu_{\mathrm{hyp}}(z)}\leq\int\limits_{\ell_{X_{0}}/4}^{\infty}
K_{\mathbb{H}}^{(1)}(t_{0};\rho)\,\mathrm{d}N_{X_{0}}(\rho;z)\,.
\end{align*}
With the notation of Lemma~4.6 of \cite{JK1}, we put $u:=\rho$, $a:=\ell_{X_{0}}/4$, and further
\begin{align*}
F(u)&:=K_{\mathbb{H}}^{(1)}(t_{0};\rho), \\[2mm]
g_{1}(u)&:=N_{X_{0}}(\rho;z), \\
g_{2}(u)&:=\frac{\sinh^{2}\big((\rho+2r)/2\big)-\sinh^{2}\big((\rho_{0}-2r)/2\big)}{\sinh^{2}(r/2)}+N_
{X_{0}}(\rho_{0};z),
\end{align*}
where $r:=\ell_{X_{0}}/4$ and $\rho_{0}:=3\ell_{X_{0}}/4$. By the latter choices for $r$ and $\rho_
{0}$, the inequalities
\begin{align*}
2r<\ell_{X_{0}}\,,\quad 2r<\rho_{0}<\ell_{X_{0}}
\end{align*}
hold, which enables us to apply Lemma~2.3~(a) of \cite{JL} to derive the inequality
\begin{align*}
g_{1}(u)\leq g_{2}(u).
\end{align*}
This in turn allows us to apply Lemma~4.6 of \cite{JK1}, in particular taking into account that $K_
{\mathbb{H}}^{(1)}(t_{0};\rho)$ is strictly monotone decreasing in $\rho$ by Proposition~A.2, namely
the inequality of Stieltjes integrals
\begin{align}
\label{16}
\int\limits_{a}^{\infty}F(u)\,\mathrm{d}g_{1}(u)+F(a)\,g_{1}(a)\leq\int\limits_{a}^{\infty}F(u)\,\mathrm
{d}g_{2}(u)+F(a)\,g_{2}(a).
\end{align}
Using the above notation, we get
\begin{align*}
&F(a)\,g_{1}(a)=K_{\mathbb{H}}^{(1)}(t_{0};\ell_{X_{0}}/4)\,N_{X_{0}}(\ell_{X_{0}}/4;z)=K_{\mathbb
{H}}^{(1)}(t_{0};\ell_{X_{0}}/4)\,,\\
&F(a)\,g_{2}(a)=K_{\mathbb{H}}^{(1)}(t_{0};\ell_{X_{0}}/4)\,\frac{\sinh^{2}(3\ell_{X_{0}}/8)-\sinh^{2}
(\ell_{X_{0}}/8)}{\sinh^{2}(\ell_{X_{0}}/8)}+K_{\mathbb{H}}^{(1)}(t_{0};\ell_{X_{0}}/4)\,.
\end{align*}
Furthermore, we compute
\begin{align*}
g_{2}(u)&=\frac{\sinh^{2}(\rho/2+\ell_{X_{0}}/4)-\sinh^{2}(\ell_{X_{0}}/8)}{\sinh^{2}(\ell_{X_{0}}/8)}+1 \\
&=\frac{\frac{1}{4}\big(e^{\rho+\ell_{X_{0}}/2}-2+e^{-\rho-\ell_{X_{0}}/2}\big)-\sinh^{2}(\ell_{X_{0}}/
8)}{\sinh^{2}(\ell_{X_{0}}/8)}+1,
\end{align*}
hence
\begin{align*}
\frac{\mathrm{d}g_{2}(u)}{\mathrm{d}u}&=\frac{\mathrm{d}}{\mathrm{d}\rho}\frac{\frac{1}{4}\big(e^
{\rho+\ell_{X_{0}}/2}-2+e^{-\rho-\ell_{X_{0}}/2}\big)-\sinh^{2}(\ell_{X_{0}}/8)}{\sinh^{2}(\ell_{X_{0}}/
8)} \\
&=\frac{\sinh(\rho+\ell_{X_{0}}/2)}{2\,\sinh^{2}(\ell_{X_{0}}/8)}\,.
\end{align*}
Inserting all of the above into \eqref{16}, we arrive at the bound
\begin{align*}
\int\limits_{\ell_{X_{0}}/4}^{\infty}K_{\mathbb{H}}^{(1)}(t_{0};\rho)\,\mathrm{d}N_{X_{0}}(\rho;z)&\leq
\int\limits_{\ell_{X_{0}}/4}^{\infty}K_{\mathbb{H}}^{(1)}(t_{0};\rho)\,\frac{\sinh(\rho+\ell_{X_{0}}/2)}{2
\,\sinh^{2}(\ell_{X_{0}}/8)}\,\mathrm{d}\rho \\
&+K_{\mathbb{H}}^{(1)}(t_{0};\ell_{X_{0}}/4)\bigg(\frac{\sinh^{2}(3\ell_{X_{0}}/8)}{\sinh^{2}(\ell_
{X_{0}}/8)}-1\bigg).
\end{align*}
Observing the inequality
\begin{align*}
\frac{\mu_{\mathrm{can}}(z)}{\mu_{\mathrm{shyp}}(z)}\leq 4\pi\,\frac{g_{X}\mu_{\mathrm{can}}(z)}
{\mu_{\mathrm{hyp}}(z)}
\end{align*}
proves the claimed bound.
\end{proof}
\end{nn}

\begin{nn}\label{4.4} 
\textbf{Proposition.} \emph{Let $X\longrightarrow X_{0}$ be an unramified covering of finite degree
with $X_{0}:=\Gamma_{0}\backslash\mathbb{H}$ a compact Riemann surface of genus $g_{0}>1$,
and let $\ell_{X_{0}}$ denote the length of the shortest closed geodesic on $X_{0}$. Then, the quantity
$S_{X}$ can be bounded as
\begin{align*}
S_{X}\leq\frac{D_{2}\,e^{\ell_{X_{0}}/2}}{(1-e^{-\ell_{X_{0}}/4})^{5/2}}
\end{align*}
with an absolute constant $D_{2}>0$, which can be taken to be $1.2\cdot 10^{3}$.}
\begin{proof}
We work from the estimate \eqref{15} for $S_{X}$ given in Lemma~\ref{4.3} and insert therein the
bound \eqref{10} for $K_{\mathbb{H}}^{(1)}(t_{0};\rho)$ obtained in Lemma \ref{4.2}, which we rewrite
as
\begin{align*}
K_{\mathbb{H}}^{(1)}(t_{0};\rho)\leq A_{1}(t_{0};\rho)+A_{2}(t_{0};\rho)+A_{3}(t_{0};\rho),
\end{align*}
where
\begin{align*}
A_{1}(t_{0};\rho)&:=\frac{17\sqrt{2}e^{-t_{0}/4}}{(4\pi t_{0})^{3/2}}\frac{(\rho+\log(4))e^{-\rho^{2}/
(4t_{0})}}{\sinh^{1/2}(\rho)}\,, \\
A_{2}(t_{0};\rho)&:=\frac{4\,\sqrt{2}\,e^{-(\rho/(2\sqrt{t_{0}})+\sqrt{t_{0}}/2)^{2}}}{\pi^{3/2}\,\sqrt{t_
{0}}}\,, \\[1mm]
A_{3}(t_{0};\rho)&:=\frac{4\,\sqrt{2}\,e^{-\rho}}{\pi^{3/2}}\int\limits_{\rho/(2\sqrt{t_{0}})-\sqrt{t_{0}}/2}^
{\infty}e^{-r'^{2}}\,\mathrm{d}r'\,.
\end{align*}
With this notation and keeping in mind that our bounds are valid for all $t_{0}>0$, we can rewrite
\eqref{15} in the form
\begin{align*}
S_{X}\leq B_{1}(t_{0};\ell_{X_{0}})+B_{2}(t_{0};\ell_{X_{0}})+B_{3}(t_{0};\ell_{X_{0}}),
\end{align*}
where
\begin{align*}
&B_{j}(t_{0};\ell_{X_{0}}):=4\pi\,\int\limits_{\ell_{X_{0}}/4}^{\infty}A_{j}(t_{0};\rho)\,\frac{\sinh(\rho+
\ell_{X_{0}}/2)}{2\,\sinh^{2}(\ell_{X_{0}}/8)}\,\mathrm{d}\rho+4\pi\,A_{j}(t_{0};\ell_{X_{0}}/4)\bigg
(\frac{\sinh^{2}(3\ell_{X_{0}}/8)}{\sinh^{2}(\ell_{X_{0}}/8)}-1\bigg)
\end{align*}
for $j=1,\,2,\,3$. In order to obtain a precise, effective upper bound for $S_{X}$, we will evaluate
the expression under consideration at $t_{0}=10$; there is no particular reason for this choice
of $t_{0}$ except to derive an explicit bound for $S_{X}$.

For the first summand  of $B_{1}(t_{0};\ell_{X_{0}})$ involving the integral, since $\sinh(\rho+
\ell_{X_{0}}/2)\leq e^{\rho+\ell_{X_{0}}/2}$ and
\begin{align*}
\frac{1}{\sinh^{1/2}(\rho)}=\frac{\sqrt{2}\,e^{-\rho/2}}{(1-e^{-2\rho})^{1/2}}\leq\frac{\sqrt{2}\,e^{-
\rho/2}}{(1-e^{-\ell_{X_{0}}/2})^{1/2}}
\end{align*}
for $\rho\geq\ell_{X_{0}}/4$, we have the bound
\begin{align*}
&\int\limits_{\ell_{X_{0}}/4}^{\infty}A_{1}(t_{0};\rho)\,\frac{\sinh(\rho+\ell_{X_{0}}/2)}{2\,\sinh^{2}
(\ell_{X_{0}}/8)}\,\mathrm{d}\rho\leq \\
&\frac{17\,e^{\ell_{X_{0}}/2}}{(4\pi t_{0})^{3/2}\sinh^{2}(\ell_{X_{0}}/8)(1-e^{-\ell_{X_{0}}/2})^{1/2}}
\int\limits_{\ell_{X_{0}}/4}^{\infty}(\rho+\log(4))e^{-(\rho/(2\sqrt{t_{0}})-\sqrt{t_{0}}/2)^{2}}\mathrm
{d}\rho\leq \\
&\frac{34\,e^{\ell_{X_{0}}/2}}{(4\pi t_{0})^{3/2}\sinh^{2}(\ell_{X_{0}}/8)(1-e^{-\ell_{X_{0}}/2})^{1/2}}
\int\limits_{-\infty}^{\infty}\big(t_{0}^{3/2}+2t_{0}\vert\rho'\vert+\log(4)t_{0}^{1/2}\big)e^{-\rho'^{2}}\,
\mathrm{d}\rho'= \\[2mm]
&\frac{34\,e^{\ell_{X_{0}}/2}}{(4\pi)^{3/2}\sinh^{2}(\ell_{X_{0}}/8)(1-e^{-\ell_{X_{0}}/2})^{1/2}}
\bigg(\bigg(1+\frac{\log(4)}{t_{0}}\bigg)\sqrt{\pi}+\frac{2}{\sqrt{t_{0}}}\bigg),
\end{align*}
hence we obtain for $t_{0}=10$
\begin{align*}
&\int\limits_{\ell_{X_{0}}/4}^{\infty}A_{1}(10;\rho)\,\frac{\sinh(\rho+\ell_{X_{0}}/2)}{2\,\sinh^{2}
(\ell_{X_{0}}/8)}\,\mathrm{d}\rho\leq\frac{3\,e^{\ell_{X_{0}}/2}}{\sinh^{2}(\ell_{X_{0}}/8)(1-e^
{-\ell_{X_{0}}/2})^{1/2}}\,.
\end{align*}
We thus get the bound
\begin{align}
\notag
\frac{B_{1}(10,\ell_{X_{0}})}{4\pi}&\leq\frac{3\,e^{\ell_{X_{0}}/2}}{\sinh^{2}(\ell_{X_{0}}/8)(1-e^
{-\ell_{X_{0}}/2})^{1/2}}+A_{1}(10;\ell_{X_{0}}/4)\bigg(\frac{\sinh^{2}(3\ell_{X_{0}}/8)}{\sinh^{2}
(\ell_{X_{0}}/8)}-1\bigg) \\
\notag
&\leq\frac{3\,e^{\ell_{X_{0}}/2}}{\sinh^{2}(\ell_{X_{0}}/8)(1-e^{-\ell_{X_{0}}/2})^{1/2}}+\frac{17\,
\sqrt{2}\,(1+\log(4))}{(40\,\pi)^{3/2}\sinh^{1/2}(\ell_{X_{0}}/4)}\frac{\sinh^{2}(3\ell_{X_{0}}/8)}{\sinh^
{2}(\ell_{X_{0}}/8)} \\
\label{17}
&\leq\frac{3\,e^{\ell_{X_{0}}/2}}{\sinh^{2}(\ell_{X_{0}}/8)(1-e^{-\ell_{X_{0}}/2})^{1/2}}+\frac{e^{5
\ell_{X_{0}}/8}}{\sinh^{2}(\ell_{X_{0}}/8)(1-e^{-\ell_{X_{0}}/2})^{1/2}}\,.
\end{align}
For the first summand of $B_{2}(t_{0};\ell_{X_{0}})$ involving the integral, we have the bound
\begin{align*}
\int\limits_{\ell_{X_{0}}/4}^{\infty}A_{2}(t_{0};\rho)\,\frac{\sinh(\rho+\ell_{X_{0}}/2)}{2\,\sinh^{2}
(\ell_{X_{0}}/8)}\,\mathrm{d}\rho\leq\frac{2\,\sqrt{2}\,e^{\ell_{X_{0}}/2}}{\pi^{3/2}\,\sqrt{t_{0}}\,
\sinh^{2}(\ell_{X_{0}}/8)}\int\limits_{\ell_{X_{0}}/4}^{\infty}e^{-(\rho/(2\sqrt{t_{0}})-\sqrt{t_{0}}/2)^
{2}}\,\mathrm{d}\rho\,,
\end{align*}
hence we obtain for $t_{0}=10$
\begin{align*}
\int\limits_{\ell_{X_{0}}/4}^{\infty}A_{2}(10;\rho)\,\frac{\sinh(\rho+\ell_{X_{0}}/2)}{2\,\sinh^{2}
(\ell_{X_{0}}/8)}\,\mathrm{d}\rho\leq\frac{4\,\sqrt{2}\,e^{\ell_{X_{0}}/2}}{\pi\,\sinh^{2}(\ell_{X_
{0}}/8)}\leq\frac{2\,e^{\ell_{X_{0}}/2}}{\sinh^{2}(\ell_{X_{0}}/8)}\,.
\end{align*}
We thus get the bound
\begin{align}
\notag
\frac{B_{2}(10,\ell_{X_{0}})}{4\pi}&\leq\frac{2\,e^{\ell_{X_{0}}/2}}{\sinh^{2}(\ell_{X_{0}}/8)}+A_
{2}(10;\ell_{X_{0}}/4)\bigg(\frac{\sinh^{2}(3\ell_{X_{0}}/8)}{\sinh^{2}(\ell_{X_{0}}/8)}-1\bigg) \\
\label{18}
&\leq\frac{2\,e^{\ell_{X_{0}}/2}}{\sinh^{2}(\ell_{X_{0}}/8)}+\frac{\sinh^{2}(3\ell_{X_{0}}/8)}{\sinh^
{2}(\ell_{X_{0}}/8)}\leq\frac{3\,e^{3\ell_{X_{0}}/4}}{\sinh^{2}(\ell_{X_{0}}/8)}\,.
\end{align}
For the first summand  of $B_{3}(t_{0};\ell_{X_{0}})$ involving the integral, we have the bound
\begin{align*}
&\int\limits_{\ell_{X_{0}}/4}^{\infty}A_{3}(t_{0};\rho)\,\frac{\sinh(\rho+\ell_{X_{0}}/2)}{2\,\sinh^{2}
(\ell_{X_{0}}/8)}\,\mathrm{d}\rho\leq \\
&\frac{2\,\sqrt{2}\,e^{\ell_{X_{0}}/2}}{\pi^{3/2}\,\sinh^{2}(\ell_{X_{0}}/8)}\int\limits_{\ell_{X_{0}}/
4}^{\infty}\int\limits_{\rho/(2\sqrt{t_{0}})-\sqrt{t_{0}}/2}^{\infty}e^{-r'^{2}}\,\mathrm{d}r'\,\mathrm
{d}\rho= \\
&\frac{2\,\sqrt{2}\,e^{\ell_{X_{0}}/2}}{\pi^{3/2}\,\sinh^{2}(\ell_{X_{0}}/8)}\int\limits_{\ell_{X_{0}}/
(8\sqrt{t_{0}})-\sqrt{t_{0}}/2}^{\infty}\int\limits_{\ell_{X_{0}}/4}^{2\sqrt{t_{0}}r'+t_{0}}e^{-r'^{2}}\,
\mathrm{d}\rho\,\mathrm{d}r'= \\
&\frac{2\,\sqrt{2}\,e^{\ell_{X_{0}}/2}}{\pi^{3/2}\,\sinh^{2}(\ell_{X_{0}}/8)}\int\limits_{\ell_{X_{0}}/
(8\sqrt{t_{0}})-\sqrt{t_{0}}/2}^{\infty}\bigg(2\sqrt{t_{0}}r'+t_{0}-\frac{\ell_{X_{0}}}{4}\bigg)e^{-r'^
{2}}\,\mathrm{d}r'\leq \\
&\frac{2\,\sqrt{2}\,e^{\ell_{X_{0}}/2}}{\pi^{3/2}\,\sinh^{2}(\ell_{X_{0}}/8)}\big(2\,\sqrt{t_{0}}+\sqrt
{\pi}\,t_{0}\big)\,,
\end{align*}
hence we obtain for $t_{0}=10$
\begin{align*}
&\int\limits_{\ell_{X_{0}}/4}^{\infty}A_{3}(10;\rho)\,\frac{\sinh(\rho+\ell_{X_{0}}/2)}{2\,\sinh^{2}
(\ell_{X_{0}}/8)}\,\mathrm{d}\rho\leq\frac{13\,e^{\ell_{X_{0}}/2}}{\sinh^{2}(\ell_{X_{0}}/8)}\,.
\end{align*}
We thus get the bound
\begin{align}
\notag
\frac{B_{3}(10,\ell_{X_{0}})}{4\pi}&\leq\frac{13\,e^{\ell_{X_{0}}/2}}{\sinh^{2}(\ell_{X_{0}}/8)}+
A_{3}(10;\ell_{X_{0}}/4)\bigg(\frac{\sinh^{2}(3\ell_{X_{0}}/8)}{\sinh^{2}(\ell_{X_{0}}/8)}-1\bigg) \\
\label{19}
&\leq\frac{13\,e^{\ell_{X_{0}}/2}}{\sinh^{2}(\ell_{X_{0}}/8)}+\frac{2\,\sinh^{2}(3\ell_{X_{0}}/
8)}{\sinh^{2}(\ell_{X_{0}}/8)}\leq\frac{15\,e^{3\ell_{X_{0}}/4}}{\sinh^{2}(\ell_{X_{0}}/8)}\,.
\end{align}
Adding up the bounds \eqref{17} -- \eqref{19}, we obtain
\begin{align*}
S_{X}&\leq\frac{88\,\pi\,e^{3\ell_{X_{0}}/4}}{\sinh^{2}(\ell_{X_{0}}/8)(1-e^{-\ell_{X_{0}}/2})^
{1/2}} \\
&\leq\frac{352\,\pi\,e^{\ell_{X_{0}}/2}}{(1-e^{-\ell_{X_{0}}/4})^{5/2}}\,,
\end{align*}
which proves the claim.
\end{proof}
\end{nn}

\begin{nn}\label{4.5} 
\textbf{Remark.} In addition to the cartesian coordinates $x$, $y$, we introduce the euclidean
polar coordinates $\rho=\rho(z)$, $\theta=\theta(z)$ of the point $z$ centered at the origin.
These are related to $x$, $y$ by the formulae
\begin{align}
\label{20}
x:=e^{\rho}\cos(\theta)\,,\quad y:=e^{\rho}\sin(\theta).
\end{align}
Given $\gamma\in\mathcal{H}(\Gamma)$, then there exists $\sigma_{\gamma}\in\mathrm{PSL}_
{2}(\mathbb{R})$ such that
\begin{align*}
\sigma_{\gamma}^{-1}\gamma\,\sigma_{\gamma}=\bigg(\begin{matrix}e^{\ell_{\gamma}/2}&
0\\0&e^{-\ell_{\gamma}/2}\end{matrix}\bigg)\quad\Longleftrightarrow\quad\gamma=\sigma_
{\gamma}\bigg(\begin{matrix}e^{\ell_{\gamma}/2}&0\\0&e^{-\ell_{\gamma}/2}\end{matrix}\bigg)
\sigma_{\gamma}^{-1}\,.
\end{align*}
For $s\in\mathbb{C}$, $\mathrm{Re}(s)>1$, the hyperbolic Eisenstein series $\mathcal{E}_
{\mathrm{hyp},\gamma}(z,s)$ associated to $\gamma$ is defined by the series
\begin{align}
\label{21}
\mathcal{E}_{\mathrm{hyp},\gamma}(z,s):=\sum\limits_{\eta\in\langle\gamma\rangle\backslash
\Gamma}\sin\big(\theta(\sigma_{\gamma}^{-1}\eta z)\big)^{s}
\end{align}
using the polar coordinates \eqref{20}. The hyperbolic Eisenstein series \eqref{21} is absolutely
and locally uniformly convergent for $z\in\mathbb{H}$ and $s\in\mathbb{C}$ with $\mathrm{Re}
(s)>1$; it is invariant under the action of $\Gamma$ and satisfies the differential equation
\begin{align*}
\big(\Delta_{\mathrm{hyp}}-s(1-s)\big)\mathcal{E}_{\mathrm{hyp},\gamma}(z,s)=s^{2}\,\mathcal
{E}_{\mathrm{hyp},\gamma}(z,s+2).
\end{align*}
For proofs of theses facts and further details, we refer to \cite{KM}.

By means of the hyperbolic Eisenstein series the following alternative bound for the quantity
$S_{X}$, namely
\begin{align*}
S_{X}\leq 8\bigg[\sum\limits_{\gamma\in\mathcal{H}(\Gamma_{0})}\bigg(\sup\limits_{z\in X_
{0}}\big\vert\Delta_{\mathrm{hyp}}\mathcal{E}_{\mathrm{hyp},\gamma}(z,2)\big\vert\,e^{-\ell_
{\gamma}}+\frac{24}{(1-e^{-\ell_{X_{0}}})^{3}}\sup\limits_{z\in X_{0}}\big\vert\mathcal{E}_
{\mathrm{hyp},\gamma}(z,2)\big\vert\,e^{-2\ell_{\gamma}}\bigg)+170\bigg],
\end{align*}
has been obtained in \cite{JK5}. This upper bound for $S_{X}$ involves special values of
hyperbolic Eisenstein series in the half-plane of convergence of the series. As such, it is
possible to use various counting function arguments, as above, to complete this approach
to obtaining an upper bound for the quantity $S_{X}$ analogous to the one given in
Proposition~\ref{4.4}.
\end{nn}


\section{Effective bounds for the Huber constant}

\begin{nn}\label{5.1} 
\textbf{Remark.} In Table $2$ of the recent joint work \cite{FJK} with J.~S.~Friedman, an algorithm
was given to bound the Huber constant $C_{\mathrm{Hub},X}$ for $X$ effectively in terms of
our basic quantities $g_{X}$, $d_{X}$, $\ell_{X}$, $\lambda_{X,1}$, and $N_{\mathrm{ev},X}^
{[0,1/4)}$; here the newly introduced quantity $d_{X}$ denotes the diameter of $X$. In the
subsequent proposition, we will summarize the result of this algorithm by utilizing convenient
yet possibly crude estimates.
\end{nn}

\begin{nn}\label{5.2} 
\textbf{Proposition.}  \emph{The Huber constant $C_{\mathrm{Hub},X}$ for $X$ can be bounded
as
\begin{align*}
C_{\mathrm{Hub},X}\leq\frac{D_{3}\,g_{X}\,e^{8\pi g_{X}/\ell_{X}+\ell_{X}/2}}{(1-s_{X,1})(1-e^
{-\ell_{X}/2})^{2}}\,;
\end{align*}
here $\ell_{X}$ denotes the length of the shortest closed geodesic on $X$,
\begin{align*}
s_{X,1}:=\frac{1}{2}+\sqrt{\frac{1}{4}-\lambda_{X,1}}
\end{align*}
with $\lambda_{X,1}$ denoting the smallest non-zero eigenvalue of $\Delta_{\mathrm{hyp}}$,
and $D_{3}>0$ is an absolute constant, which can be taken to be $10^{11}$.}
\begin{proof}
As mentioned in \ref{5.1}, we follow the algorithm given in Table $2$ of \cite{FJK}. In the sequel
we also use the definitions of the quantities $A$, $B$, $C$, $C_{j}$ ($j=1,2,3,\ldots$) therein.

Recalling  from \cite{Buser} the bound for $N_{\mathrm{ev},X}^{[0,1/4)}$, we obtain for the
quantity $A$ the estimate
\begin{align*}
A:=N_{\mathrm{ev},X}^{[0,1/4)}\leq 4g_{X}-2\leq 4g_{X}.
\end{align*}
Using the inequality~(2) from the main theorem of \cite{Chang}, namely
\begin{align*}
2\,\frac{\ell_{X}}{4}\,d_{X}\leq2\sinh\bigg(\frac{\ell_{X}}{4}\bigg)d_{X}\leq 4\pi(g_{X}-1),
\end{align*}
we obtain the following bound for the diameter $d_{X}$ of $X$
\begin{align*}
d_{X}\leq\frac{8\,\pi\,g_{X}}{\ell_{X}}.
\end{align*}
Hence, the quantity $B$ can be estimated by
\begin{align*}
B:=\frac{2\pi e^{d_{X}}}{4\pi(g_{X}-1)}\leq\frac{e^{8\pi g_{X}/\ell_{X}}}{2}\,.
\end{align*}
For the quantity $C$, we have
\begin{align*}
C:=3\bigg(\frac{4\pi(g_{X}-1)}{4\pi}+745B\bigg)\leq 3\,g_{X}+1118\,e^{8\pi g_{X}/\ell_{X}}.
\end{align*}
Next, we have
\begin{align*}
C_{1}:=2e-2\leq 4,
\end{align*}
and
\begin{align*}
C_{10}:= 8480\sqrt{\frac{e}{2\pi}}\leq 5578.
\end{align*}
From this we derive
\begin{align*}
C_{12}&:=(A-1)\bigg(1+3C_{1}+\frac{2}{1-s_{X,1}}(1+C_{1})\bigg)+2C_{1}+2 \\
&\leq 4\,g_{X}\bigg(13+\frac{10}{1-s_{X,1}}\bigg)+10 \\
&\leq\frac{92\,g_{X}}{1-s_{X,1}}+10\leq\frac{102\,g_{X}}{1-s_{X,1}}\,,
\end{align*}
and
\begin{align*}
C_{13}&:=\frac{41}{6}\cdot C\cdot C_{10}\leq\frac{41\cdot 5578}{6}\big(3\,g_{X}+1118\,e^
{8\pi g_{X}/\ell_{X}}\big) \\
&\leq 114\,349\,g_{X}+42\,614\,061\,e^{8\pi g_{X}/\ell_{X}}.
\end{align*}
From this we obtain
\begin{align*}
C_{16}&:=C_{12}+C_{13}+\frac{3}{2\pi}4\pi(g_{X}-1)C_{10} \\
&\leq\frac{102\,g_{X}}{1-s_{X,1}}+114\,349\,g_{X}+42\,614\,061\,e^{8\pi g_{X}/\ell_{X}}+6
\cdot 5578\,g_{X} \\
&\leq\frac{147\,919\,g_{X}}{1-s_{X,1}}+42\,614\,061\,e^{8\pi g_{X}/\ell_{X}}\leq\frac{42\,
761\,980\,g_{X}\,e^{8\pi g_{X}/\ell_{X}}}{1-s_{X,1}}\,.
\end{align*}
For notational convenience, let us keep the constant $C_{16}$ without replacing it with the
above bound for the next few computations. We further have
\begin{align*}
&C_{17}:=4A+4\,C_{16}\leq 16\,g_{X}+4\,C_{16}, \\
&C_{18}:=4A+5\,C_{16}\leq 16\,g_{X}+5\,C_{16}.
\end{align*}
The constant $c$ must satisfy $1<c<e^{\ell_{X}}$, so we may take $c:=e^{\ell_{X}/2}$, and
hence $\mu:=\ell_{X}/2$. With this choice, we find
\begin{align*}
C_{19}&:=C_{18}+\frac{8A+4\,C_{18}}{1-1/c}\leq 16\,g_{X}+5\,C_{16}+\frac{96\,g_{X}+20\,
C_{16}}{1-e^{-\ell_{X}/2}} \\
&\leq\frac{112\,g_{X}+25\,C_{16}}{1-e^{-\ell_{X}/2}}\,.
\end{align*}
Observing that
\begin{align*}
f(r):=\frac{r}{1-e^{-r}}\geq 1
\end{align*}
for $r\in\mathbb{R}_{\geq 0}$, we find
\begin{align*}
\frac{1}{\mu}=\frac{2}{\ell_{X}} \leq\frac{1}{1-e^{-\ell_{X}/2}}\,.
\end{align*}
Thus, we obtain
\begin{align*}
C_{20}&:=C_{19}+\frac{8A+4\,C_{18}}{\mu}\leq\frac{112\,g_{X}+25\,C_{16}}{1-e^{-\ell_{X}/
2}}+\frac{8A+4\,C_{18}}{1-e^{-\ell_{X}/2}} \\
&\leq\frac{112\,g_{X}+25\,C_{16}}{1-e^{-\ell_{X}/2}}+\frac{96\,g_{X}+20\,C_{16}}{1-e^{-\ell_
{X}/2}} \\
&=\frac{208\,g_{X}+45\,C_{16}}{1-e^{-\ell_{X}/2}}\,.
\end{align*}
For the quantity $C_{21}$, we find the estimate
\begin{align*}
C_{21}&:=\frac{\vert c-2\vert}{\log(2)}+\frac{2\vert 2-\sqrt{c}\vert}{\log(c)}\leq\frac{c+2}{\log(2)}+
\frac{2(\sqrt{c}+2)}{\log(c)} \\
&\leq\frac{e^{\ell_{X}/2}+2}{\log(2)}+\frac{4(e^{\ell_{X}/4}+2)}{\ell_{X}}\leq\frac{3\,e^{\ell_{X}/
2}}{1/2}+\frac{4\cdot 3\,e^{\ell_{X}/2}}{\ell_{X}} \\
&\leq 12\,e^{\ell_{X}/2}\bigg(1+\frac{1}{\ell_{X}}\bigg)\leq\frac{18\,e^{\ell_{X}/2}}{1-e^{-\ell_
{X}/2}}\,.
\end{align*}
At this point, we have to correct the statement about the constant $C_{22}$, which comes
from Lemma~4.14 in \cite{FJK}. The correct assertion is that
\begin{align*}
C_{22}:=\frac{1}{1+1/\log(2)}\,.
\end{align*}
In fact, $C_{22}$ has to be such that for any $r\geq 2$, we have the inequality
\begin{align*}
\mathrm{li}(r)\leq C_{22}\,\frac{r}{\log(r)}\,.
\end{align*}
For a proof we consider the function
\begin{align*}
f(r):=\mathrm{li}(r)-d\,\frac{r}{\log(r)}
\end{align*}
for some positive constant $d$, which we determine such that $f(r)$ is negative for $r\geq 2$.
Obviously, $f(2)<0$, so we have to determine $d$ such that $f(r)$ becomes a decreasing
function. We have
\begin{align*}
f'(r)=\frac{1}{\log(r)}\bigg(1-d+\frac{d}{\log(r)}\bigg),
\end{align*}
hence, we need to have
\begin{align*}
1-d+\frac{d}{\log(r)}\leq 0\quad\Longleftrightarrow\quad 1-\frac{1}{d}\geq\frac{1}{\log(r)}
\end{align*}
for $r\geq 2$, which holds for
\begin{align*}
1-\frac{1}{d}\geq\frac{1}{\log(2)}\quad\Longleftrightarrow\quad d\geq\frac{1}{1+1/\log(2)}
\end{align*}
giving the claimed value of $C_{22}$. (Note that the error in the proof of Lemma~4.14 of \cite
{FJK} arose by dividing by a constant which is negative, so then the inequality has to change
directions.) Continuing with this value of $C_{22}$, we have
\begin{align*}
C_{22}=\frac{1}{1+1/\log(2)}\leq\frac{1}{2}.
\end{align*}
Finally, we are in a position to compute $C_{u}$; we have
\begin{align*}
C_{u}&:=C_{21}A+C_{20}\frac{c^{3/4}}{\log(c)}+C_{20}(1+C_{22})+\frac{3}{4}C_{20}C_{21} \\
&\leq\frac{72\,g_{X}\,e^{\ell_{X}/2}}{1-e^{-\ell_{X}/2}}+\frac{(208\,g_{X}+45\,C_{16})e^{\ell_{X}/
2}}{(1-e^{-\ell_{X}/2})^{2}} \\
&+\frac{312\,g_{X}+69\,C_{16}}{1-e^{-\ell_{X}/2}}+\frac{(3744\,g_{X}+810\,C_{16})e^{\ell_{X}/
2}}{(1-e^{-\ell_{X}/2})^{2}}\,.
\end{align*}
Employing finally the bound for $C_{16}$ yields the estimate
\begin{align*}
C_{u}&\leq\frac{384\,g_{X}\,e^{\ell_{X}/2}}{1-e^{-\ell_{X}/2}}+\frac{69\,C_{16}}{1-e^{-\ell_{X}/
2}}+\frac{(3952\,g_{X}+855\,C_{16})e^{\ell_{X}/2}}{(1-e^{-\ell_{X}/2})^{2}} \\
&\leq\frac{39\,512\,073\,856\,g_{X}\,e^{8\pi g_{X}/\ell_{X}+\ell_{X}/2}}{(1-s_{X,1})(1-e^{-\ell_
{X}/2})^{2}}\,.
\end{align*}
This completes the proof of the proposition.
\end{proof}
\end{nn}

\section{Effective bounds for Faltings's delta function}

The main result proven in this paper consists in simplifying the bound obtained in Corollary~\ref
{3.4} and making it effective.

\begin{nn}\label{6.1} 
\textbf{Theorem.} \emph{Let $X\longrightarrow X_{0}$ be an unramified covering of finite degree with
$X_{0}:=\Gamma_{0}\backslash\mathbb{H}$ a compact Riemann surface of genus $g_{X_{0}}>1$.
Let $\ell_{X_{0}}$ denote the length of the shortest closed geodesic on $X_{0}$ and $\lambda_{X,1}$,
$\lambda_{X_{0},1}$ the smallest non-zero eigenvalues of $\Delta_{\mathrm{hyp}}$ on $X$, $X_{0}$,
respectively, and
\begin{align*}
\lambda_{X}=\frac{1}{2}\min\bigg\{\lambda_{X,1},\frac{7}{64}\bigg\}\,,\qquad s_{X_{0},1}=\frac{1}
{2}+\sqrt{\frac{1}{4}-\lambda_{X_{0},1}}\,.
\end{align*}
Then, we have the effective bound
\begin{align}
\label{22}
\delta_{\mathrm{Fal}}(X)\leq\frac{D_{4}\,g_{X_{0}}\,e^{8\pi g_{X_{0}}/\ell_{X_{0}}+\ell_{X_{0}}}}
{(1-e^{-\ell_{X_{0}}/4})^{5}(1-s_{X_{0},1})}\,\frac{g_{X}}{\lambda_{X}}
\end{align}
with an absolute constant $D_{4}>0$, which can be taken to be $10^{15}$.}
\begin{proof}
We work from the bound
\begin{align}
\label{23}
\delta_{\mathrm{Fal}}(X)\leq 876\bigg(g_{X}+\frac{1}{\lambda_{X}}\big(g_{X}(S_{X}+1)^{2}+C_
{\mathrm{Hub},X}+N_{\mathrm{ev},X}^{[0,1/4)}\big)+\bigg(1+\frac{1}{\ell_{X}}\bigg)N_{\mathrm
{geo},X}^{(0,5)}\bigg).
\end{align}
obtained in the proof of Corollary \ref{3.4} (using the notation therein). We will next bound the
quantities
\begin{align*}
\ell_{X},\quad N_{\mathrm{ev},X}^{[0,1/4)},\quad S_{X},\quad C_{\mathrm{Hub},X},\quad N_
{\mathrm{geo},X}^{(0,5)}
\end{align*}
in terms of the underlying compact Riemann surface $X_{0}$.

(i) We start by observing that the trivial inequality
\begin{align}
\label{24}
\ell_{X}\geq\ell_{X_{0}}
\end{align}
holds true for the lengths of the shortest closed geodesics on $X$, $X_{0}$, respectively.

(ii) In order to estimate $N_{\mathrm{ev},X}^{[0,1/4)}$, we recall as in the proof of Proposition
\ref{5.2} from \cite{Buser} the bound
\begin{align}
\label{25}
1\leq N_{\mathrm{ev},X}^{[0,1/4)}\leq 4g_{X}-2\leq 4g_{X}.
\end{align}

(iii) From Proposition \ref{4.4}, we recall the bound
\begin{align}
\label{26}
S_{X}\leq\frac{1\,200\,e^{\ell_{X_{0}}/2}}{(1-e^{-\ell_{X_{0}}/4})^{5/2}}
\end{align}
with $\ell_{X_{0}}$ as in the statement of the theorem.

(iv) Next, we have to estimate $C_{\mathrm{Hub},X}$. We start by citing Theorem~3.4 of \cite{JK2} 
and use the Artin formalism for the covering $X\longrightarrow X_{0}$, to derive the bound
\begin{align*}
C_{\mathrm{Hub},X}\leq[\Gamma_{0}\colon\Gamma]\,C_{\mathrm{Hub},X_{0}}.
\end{align*}
From the Riemann--Hurwitz formula we now easily derive the bound
\begin{align*}
[\Gamma_{0}\colon\Gamma]\leq\frac{g_{X}-1}{g_{X_{0}}-1}\leq g_{X},
\end{align*}
from which we get
\begin{align}
\label{27}
C_{\mathrm{Hub},X}\leq g_{X}\,C_{\mathrm{Hub},X_{0}},
\end{align}
where the proof of Proposition \ref{5.2} shows
\begin{align}
\label{28}
C_{\mathrm{Hub},X_{0}}\leq\frac{39\,512\,073\,856\,g_{X_{0}}\,e^{8\pi g_{X_{0}}/\ell_{X_{0}}+
\ell_{X_{0}}/2}}{(1-s_{X_{0},1})(1-e^{-\ell_{X_{0}}/2})^{2}}
\end{align}
with $\ell_{X_{0}}$ and $s_{X_{0},1}$ as in the statement of the theorem.

(v) Finally, we need to bound $N_{\mathrm{geo},X}^{(0,5)}$. With the above notation, using
arguments from the proof of Theorem~4.11 in \cite{JK1} (as well as the notation $r_{\Gamma_
{0},\Gamma}$ therein), we find (as above)
\begin{align*}
N_{\mathrm{geo},X}^{(0,5)}\leq\frac{5\,r_{\Gamma_{0},\Gamma}}{\ell_{X_{0}}}\,N_{\mathrm{geo},
X_{0}}^{(0,5)}\leq\frac{5\,[\Gamma_{0}\colon\Gamma]}{\ell_{X_{0}}}\,N_{\mathrm{geo},X_{0}}^
{(0,5)}
\leq\frac{5\,g_{X}}{\ell_{X_{0}}}\,N_{\mathrm{geo},X_{0}}^{(0,5)}\,.
\end{align*}
Applying the prime geodesic theorem \eqref{1} to $X_{0}$ and recalling the monotonicity of
the logarithmic integral for $u>0$, we find
\begin{align}
\notag
N_{\mathrm{geo},X_{0}}^{(0,5)}=\pi_{X_{0}}\big(\log(5)\big)&\leq N_{\mathrm{ev},X_{0}}^{[0,1/4)}\,
\mathrm{li}\big(\log(5)\big)+C_{\mathrm{Hub},X_{0}}\,\frac{\log(5)^{3/4}}{\log\big(\log(5)\big)^{1/2}}
\\
\label{29}
&\leq N_{\mathrm{ev},X_{0}}^{[0,1/4)}+3\,C_{\mathrm{Hub},X_{0}}\leq 4g_{X_{0}}+3\,C_{\mathrm
{Hub},X_{0}}\,,
\end{align}
where $C_{\mathrm{Hub},X_{0}}$ can be effectively bounded using Proposition \ref{5.2}.

Inserting the bounds \eqref{24} -- \eqref{29} into the estimate \eqref{23} yields the following
bound for $\delta_{\mathrm{Fal}}(X)$:
\begin{align*}
& 876\Bigg(g_{X}+\frac{1}{\lambda_{X}}\bigg(g_{X}\bigg(\frac{1\,200\,e^{\ell_{X_{0}}/2}}{(1-e^
{-\ell_{X_{0}}/4})^{5/2}}+1\bigg)^{2}+g_{X}C_{\mathrm{Hub},X_{0}}+4g_{X}\bigg)+\frac{5\,g_{X}}
{\ell_{X_{0}}}\bigg(1+\frac{1}{\ell_{X_{0}}}\bigg)N_{\mathrm{geo},X_{0}}^{(0,5)}\Bigg)\leq \\
&\hspace*{1cm}876\,g_{X}\bigg(1+\frac{1}{\lambda_{X}}\bigg(\frac{1\,201^{2}\,e^{\ell_{X_{0}}}}
{(1-e^{-\ell_{X_{0}}/4})^{5}}+C_{\mathrm{Hub},X_{0}}+4\bigg)+\frac{10\,N_{\mathrm{geo},X_
{0}}^{(0,5)}}{\ell_{X_{0}}(1-e^{-\ell_{X_{0}}/2})}\bigg)\leq \\
&\hspace*{1cm}876\,g_{X}\bigg(\frac{1}{\lambda_{X}}\bigg(\frac{1\,442\,401\,e^{\ell_{X_{0}}}}
{(1-e^{-\ell_{X_{0}}/4})^{5}}+C_{\mathrm{Hub},X_{0}}+5\bigg)+\frac{20\,g_{X_{0}}+15\,C_
{\mathrm{Hub},X_{0}}}{(1-e^{-\ell_{X_{0}}/2})^{2}}\bigg)\leq \\
&\hspace*{1cm}876\,\frac{g_{X}}{\lambda_{X}}\bigg(\frac{1\,442\,426\,g_{X_{0}}\,e^{\ell_{X_
{0}}}}{(1-e^{-\ell_{X_{0}}/4})^{5}}+\frac{16\,C_{\mathrm{Hub},X_{0}}}{(1-e^{-\ell_{X_{0}}/4})^
{2}}\bigg)\leq \\
&\hspace*{1cm}\frac{876}{(1-e^{-\ell_{X_{0}}/4})^{5}}\frac{g_{X}}{\lambda_{X}}\bigg(1\,442\,
426\,g_{X_{0}}\,e^{\ell_{X_{0}}}+\frac{632\,193\,181\,696\,g_{X_{0}}\,e^{8\pi g_{X_{0}}/\ell_
{X_{0}}+\ell_{X_{0}}/2}}{1-s_{X_{0},1}}\bigg)\leq \\
&\hspace*{1cm}553\,802\,490\,730\,872\,\frac{g_{X_{0}}\,e^{8\pi g_{X_{0}}/\ell_{X_{0}}+\ell_
{X_{0}}}}{(1-e^{-\ell_{X_{0}}/4})^{5}(1-s_{X_{0},1})}\,\frac{g_{X}}{\lambda_{X}}\,.
\end{align*}
This completes the proof of the theorem.
\end{proof}
\end{nn}

\begin{nn}\label{6.2} 
\textbf{Remarks.} (i) We can further refine the lower bound for $\ell_{X_{0}}$ provided that $X_
{0}$ has a model defined over some number field. In fact, by B\'elyi's theorem, we then have
$X_{0}\cong\overline{\Delta_{0}\backslash\mathbb{H}}$, where $\Delta_{0}$ is a subgroup of
finite index in $\Gamma(2)$. Therefore, we have the estimate
\begin{align*}
2\cosh(\ell_{X_{0}}/2)=\vert\mathrm{tr}(\delta_{0})\vert\geq 4,
\end{align*}
where $\delta_{0}\in\Delta_{0}$ is such that $\ell_{\delta_{0}}=\ell_{X_{0}}$; this gives $\ell_
{X_{0}}\geq 2\,\mathrm{arcosh}(2)$. The factor depending on $X_{0}$ in \eqref{22} can thus
be bounded as
\begin{align}
\label{30}
\frac{D_{4}\,g_{X_{0}}\,e^{8\pi g_{X_{0}}/\ell_{X_{0}}+\ell_{X_{0}}}}{(1-e^{-\ell_{X_{0}}/4})^{5}
(1-s_{X_{0},1})}\leq\frac{40\,D_{4}\,g_{X_{0}}\,e^{10\,g_{X_{0}}+\ell_{X_{0}}}}{1-s_{X_{0},1}}
\leq\frac{10^{17}\,g_{X_{0}}\,e^{10\,g_{X_{0}}+\ell_{X_{0}}}}{1-s_{X_{0},1}}\,.
\end{align}
(ii) On the other hand, if $X_{0}$ can be covered by a modular curve $\overline{\Gamma(N)
\backslash\mathbb{H}}$ for the full congruence subgroup $\Gamma(N)$ for some $N\in\mathbb
{N}$, a result of R.~Brooks in \cite{Brooks} shows that $\lambda_{X_{0},1}\geq 5/36$, which
gives the estimate
\begin{align*}
\frac{1}{1-s_{X_{0},1}}\leq 6.
\end{align*}
In addition, assuming that $X_{0}$ has a model defined over some number field, case (i)
above also applies and the bound \eqref{30} simplifies to
\begin{align*}
\frac{D_{4}\,g_{X_{0}}\,e^{8\pi g_{X_{0}}/\ell_{X_{0}}+\ell_{X_{0}}}}{(1-e^{-\ell_{X_{0}}/4})^{5}
(1-s_{X_{0},1})}\leq 10^{18}\,g_{X_{0}}\,e^{10\,g_{X_{0}}+\ell_{X_{0}}}\,.
\end{align*}
\end{nn}

\begin{nn}\label{6.3} 
\textbf{Corollary.} \emph{Let $X$ be a compact Riemann surface of genus $g_{X}>1$. Let $\ell_
{X}$ denote the length of the shortest closed geodesic on $X$, $\lambda_{X,1}$ the smallest
non-zero eigenvalue of $\Delta_{\mathrm{hyp}}$ on $X$, and
\begin{align*}
\lambda_{X}=\frac{1}{2}\min\bigg\{\lambda_{X,1},\frac{7}{64}\bigg\}\,,\quad s_{X,1}=\frac{1}{2}+
\sqrt{\frac{1}{4}-\lambda_{X,1}}\,.
\end{align*}
Then, we have the effective bound
\begin{align}
\label{31}
\delta_{\mathrm{Fal}}(X)\leq\frac{D_{4}\,g_{X}\,e^{8\pi g_{X}/\ell_{X}+\ell_{X}}}{(1-e^{-\ell_{X}/
4})^{5}}\frac{1}{\lambda_{X}(1-s_{X,1})}  
\end{align}
with an absolute constant $D_{4}>0$, which can be taken to be $10^{15}$.}
\begin{proof}
The proof follows immediately from an analysis of the proof of Theorem~\ref{6.1} for the trivial
covering $X_{0}=X$.
\end{proof}
\end{nn}

Using Corollary \ref{6.3}, we can now also give a variant of the bound \eqref{22} in the case
that $X$ is a \emph{ramified} covering of finite degree of a compact Riemann surface $X_{0}$
of genus $g_{X_{0}}>1$. For this, we let $\mathrm{Ram}(X/X_{0})\subset X_{0}$ denote the
ramification locus of the given covering.

\begin{nn}\label{6.4}
\textbf{Corollary.} \emph{Let $X\longrightarrow X_{0}$ be a ramified covering of finite degree
of compact Riemann surfaces of genera $g_{X},g_{X_{0}}>1$, respectively. With $\ell_{X_{0}}$
denoting the length of the shortest closed geodesic on $X_{0}$, put
\begin{align*}
r_{X_{0}}:=\min\bigg\{\ell_{X_{0}}\,,\min_{\substack{z,w\in\mathrm{Ram}(X/X_{0})\\z\neq w}}
\mathrm{dist}_{\mathrm{hyp}}(z,w)\bigg\}\,,\quad R_{X_{0}}:=\max\bigg\{\ell_{X_{0}}\,,\max_
{\substack{z,w\in\mathrm{Ram}(X/X_{0})\\z\neq w}}\mathrm{dist}_{\mathrm{hyp}}(z,w)\bigg\}\,.
\end{align*}
Furthermore, with $\lambda_{X,1}$ denoting the smallest non-zero eigenvalue of $\Delta_
{\mathrm{hyp}}$ on $X$, put
\begin{align*}
\lambda_{X}=\frac{1}{2}\min\bigg\{\lambda_{X,1},\frac{7}{64}\bigg\}\,,\quad s_{X,1}=\frac{1}
{2}+\sqrt{\frac{1}{4}-\lambda_{X,1}}\,.
\end{align*}
Then, we have the effective bound
\begin{align*}
\delta_{\mathrm{Fal}}(X)\leq\frac{D_{4}\,g_{X}\,e^{8\pi g_{X}/r_{X_{0}}+g_{X}R_{X_{0}}}}
{(1-e^{-r_{X_{0}}/4})^{5}}\frac{1}{\lambda_{X}(1-s_{X,1})}
\end{align*}
with an absolute constant $D_{4}>0$, which can be taken to be $10^{15}$.}
\begin{proof}
We work from the effective bound obtained in Corollary \ref{6.3} and estimate the length of
the shortest closed geodesic $\ell_{X}$ from below and above by quantities depending on
the base $X_{0}$.

In order to estimate $\ell_{X}$ from below, we observe that the length of closed geodesics
on $X$, which do not pass through ramification points, can be bounded from below by $\ell_
{X_{0}}$; the same estimate holds true, if the closed geodesic passes through a single 
ramification point. However, if the closed geodesic happens to pass through at least
two ramification points lying above two distinct points of $\mathrm{Ram}(X/X_{0})$, we
additionally have to take into account the distances between mutually distinct points of
$\mathrm{Ram}(X/X_{0})$ in our estimate. All in all this leads to the lower bound
\begin{align}
\label{32}
\ell_{X}\geq r_{X_{0}}.
\end{align}
Similarly, we find that the length of closed geodesics on $X$, which do not pass through
ramification points, can be bounded from above by $\deg(X/X_{0})\,\ell_{X_{0}}$, and the
same estimate holds true, if the closed geodesic passes through a single ramification
point. Again, if the closed geodesic happens to pass through at least two ramification
points lying above two distinct points of $\mathrm{Ram}(X/X_{0})$, we additionally have
to take into account the distances between mutually distinct points of $\mathrm{Ram}(X/
X_{0})$ in our estimate. This leads to the upper bound
\begin{align}
\label{33}
\ell_{X}\leq\deg(X/X_{0})\,R_{X_{0}}\leq g_{X}R_{X_{0}}.
\end{align}
Inserting the bounds \eqref{32} and \eqref{33} into the estimate \eqref{31} completes the
proof of the corollary.
\end{proof}
\end{nn}

\section{Application to Parshin's covering construction}

\begin{nn}\label{7.1} 
\textbf{The set-up.} Let $K$ be a number field with ring of integers $\mathcal{O}_{K}$ and $S:=
\mathrm{Spec}\,(\mathcal{O}_{K})$. In contrast to the previous sections, let $X$ denote a smooth
projective curve defined over $K$ of genus $g_{X}>1$, and let $\mathcal{X}/S$ be a minimal
regular model of $X/K$, which is semistable. Denote by $\overline{\omega}_{\mathcal{X}/S}$
the relative dualizing sheaf of $\mathcal{X}/S$ equipped with the Arakelov metric. For $\mathfrak
{p}\in S$, we let $\delta_{\mathfrak{p}}$ be the number of singular points in the fiber above
$\mathfrak{p}$. For an archimedean place $v$, we put
\begin{align*}
X_{v}:=X\times_{v}\mathbb{C},
\end{align*}
whose complex points $X_{v}(\mathbb{C})$ constitute a compact Riemann surface of genus
equal to $g_{X}$. In order to simplify our notation, we allow ourselves subsequently to write
$X_{v}$ instead of $X_{v}(\mathbb{C})$.

In his quest for an arithmetic version of the van de Ven--Bogomolov--Miyaoka--Yau inequality,
A.~N.~Parshin proposed the following inequality (see \cite{Parshin})
\begin{align}
\notag
&\overline{\omega}_{\mathcal{X}/S}^{2}\leq \\
\label{34}
&c_{1}\bigg(\sum\limits_{\mathfrak{p}}\delta_{\mathfrak{p}}\log\big(N_{K/\mathbb{Q}}(\mathfrak{p})
\big)+\sum\limits_{v}\varepsilon_{v}\,\delta_{\mathrm{Fal}}\big(X_{v}\big)\bigg)+c_{2}\big(2g_{X}-2
\big)\log\big\vert\mathrm{disc}(K/\mathbb{Q})\big\vert+c_{3}[K\colon\mathbb{Q}]\,;
\end{align}
here $c_{j}$ are positive constants depending solely on $K$ ($j=1,2,3$), $N_{K/\mathbb{Q}}
(\mathfrak{p})$ denotes the absolute norm of $\mathfrak{p}$, and $\mathrm{disc}(K/\mathbb{Q})$
is the discriminant of the field extension $K/\mathbb{Q}$. As is well known by subsequent work of
J.-B.~Bost, J.-F.~Mestre, and L.~Moret-Bailly (see \cite{BMMB}), the inequality \eqref{34} does not
hold true in general.
\end{nn}

\begin{nn}\label{7.2} 
\textbf{The covering construction.} Assuming the validity of the inequality \eqref{34}, A.~N.~Parshin
proposed in \cite{Parshin}, how to bound the height of $K$-rational points $P\in X(K)$ as effective
as possible using the following ramified covering construction.

Given the smooth projective curve $X/K$ of genus $g_{X}>1$, and $P\in X(K)$ a $K$-rational
point, there exists a finite covering $X_{P}/K_{P}$ over $X$ with the following properties:
\begin{itemize}
\item[(i)]
The field extension $K_{P}/K$ is a finite extension of degree effectively bounded as $O(g_{X})$
with prescribed ramification.
\item[(ii)]
The covering $X_{P}/X$ is finite of degree effectively bounded as $O(g_{X})$ and ramified only
at $P$ of ramification index effectively bounded as $O(g_{X})$; by the Riemann--Hurwitz formula,
the genus $g_{X_{P}}$ of $X_{P}$ is then also  effectively bounded as $O(g_{X})$.
\item[(iii)]
For each archimedean place $v$ of $K$ and each archimedean place $v'$ of $K_{P}$ lying
above $v$, there exists a smooth projective complex surface $Y_{v}$ together with a smooth
morphism $\varphi_{v}\colon Y_{v}\longrightarrow X_{v}$ such that
\begin{align*}
\varphi_{v}^{-1}(P)\cong X_{P,v'}:=X_{P}\times_{v'}\mathbb{C}.
\end{align*}
\end{itemize}
Denoting by $\mathcal{O}_{K_{P}}$ the ring of integers of $K_{P}$, setting $S_{P}:=\mathrm
{Spec}\,(\mathcal{O}_{K_{P}})$, letting $\mathcal{X}_{P}/S_{P}$ be a minimal regular model
of $X_{P}/K_{P}$, which is semistable, and denoting by  $\overline{\omega}_{\mathcal{X}_{P}/
S_{P}}$ the relative dualizing sheaf of $\mathcal{X}_{P}/S_{P}$ equipped with the Arakelov
metric, the height $h(P)$ of $P$ can be bounded by the arithmetic self-intersection number
\begin{align}
\label{35}
h(P)\ll\overline{\omega}_{\mathcal{X}_{P}/S_{P}}^{2}\,,
\end{align}
which, in turn, can then be bounded using \eqref{34}, after replacing $\overline{\omega}_
{\mathcal{X}/S}$ by $\overline{\omega}_{\mathcal{X}_{P}/S_{P}}$. In \cite{Parshin}, the quantities
$\delta_{\mathfrak{P}}$ ($\mathfrak{P}\in S_{P}$), $\mathrm{disc}(K_{P}/\mathbb{Q})$, and $[K_
{P}\colon\mathbb{Q}]$ are then effectively bounded in terms of the genus $g_{X}$ of $X$. The
contribution from Faltings's delta function $\delta_{\mathrm{Fal}}\big(X_{P,v'}\big)$ ($v'\vert v$)
is bounded in terms of $X$ by arguing that, as $P$ is moving through the set of $K$-rational
points $X(K)$, the function $\delta_{\mathrm{Fal}}\big(X_{P,v'}\big)$ can be viewed as the
restriction of a real-analytic function on $X_{v}$, which takes its maximum on the compact
Riemann surface $X_{v}$.
\end{nn}

\begin{nn}\label{7.3} 
\textbf{Parshin's question.} After having presented our estimate \eqref{23} for Faltings's delta
function obtained in Corollary \ref{3.4}, Parshin proposed to apply our bound to $\delta_
{\mathrm{Fal}}\big(X_{P,v'}\big)$ in order to obtain a more explicit bound than his.

Indeed, applying the bound obtained in Corollary \ref{6.4} to the ramified covering $X_{P,v'}
\longrightarrow X_{v}$ of finite degree, observing that the ramification locus $\mathrm{Ram}
(X_{P,v'}/X_{v})$ consists of only one point, we are led to the bound
\begin{align}
\label{36}
\delta_{\mathrm{Fal}}\big(X_{P,v'}\big)\leq\frac{D_{4}\,g_{X_{P}}\,e^{8\pi g_{X_{P}}/\ell_{X_
{v}}+g_{X_{P}}\ell_{X_{v}}}}{(1-e^{-\ell_{X_{v}}/4})^{5}}\frac{1}{\lambda_{X_{P,v'}}(1-s_{X_
{P,v'},1})}\,,
\end{align}
where
\begin{align*}
\lambda_{X_{P,v'}}=\frac{1}{2}\min\bigg\{\lambda_{X_{P,v'},1},\frac{7}{64}\bigg\}\,,\qquad
s_{X_{P,v'},1}=\frac{1}{2}+\sqrt{\frac{1}{4}-\lambda_{X_{P,v'},1}}
\end{align*}
with $\ell_{X_{v}}$ denoting the length of the shortest closed geodesic on $X_{v}$ and
$\lambda_{X_{P,v'},1}$ denoting the smallest non-zero eigenvalue of $\Delta_{\mathrm
{hyp}}$ on $X_{P,v'}$.

As $P$ is moving through the set of $K$-rational points $X(K)$ or, more generally, through
the compact Riemann surface $X_{v}$, the Riemann surfaces $X_{P,v'}$ (or, rather their
isomorphism classes) cover a compact region $\mathcal{D}$ in the moduli space $\mathcal
{M}_{g_{X_{P}}}$ of curves of genus $g_{X_{P}}$. While $P$ is ranging over $X_{v}$, the
function
\begin{align*}
\lambda_{X_{P,v'}}(1-s_{X_{P,v'},1})
\end{align*}
takes its minimum on $\mathcal{D}$, which we denote by $\lambda_{v,\mathrm{min}}$.
Keeping in mind that $X_{v}$ is defined over a number field, Remark~\ref{6.2}~(i) allows
us to simplify the bound \eqref{36} to
\begin{align*}
\delta_{\mathrm{Fal}}\big(X_{P,v'}\big)\leq\frac{10^{17}\,g_{X_{P}}\,e^{10g_{X_{P}}+g_
{X_{P}}\ell_{X_{v}}}}{\lambda_{v,\mathrm{min}}}\,;
\end{align*}
here we recall that the genus $g_{X_{P}}$ can be effectively bounded in terms of the
genus $g_{X}$.

We conclude by emphasizing that our results do not lead to an \emph{effective} bound
for the height $h(P)$ of $K$-rational points $P\in X(K)$, since the bound \eqref{35} as
well as the determination of the minimum $\lambda_{v,\mathrm{min}}$ are not effective.
\end{nn}

\newpage

\setcounter{secnumdepth}{0}
\section{A. Appendix}

In order to apply the inequality of Stieltjes integrals \eqref{16}, we need that the function
$K_{\mathbb{H}}^{(1)}(t;\rho)$ is monotone decreasing in $\rho$. The purpose of this
appendix is to provide a proof of this claim.

\textbf{A.1. Lemma.} \emph{For $t > 0$, $\rho>0$, and $r\geq\rho$, let
\begin{align*}
F(t;\rho,r):=\frac{re^{-r^{2}/(4t)}}{\sinh(r)}T_{2}\bigg(\frac{\cosh(r/2)}{\cosh(\rho/2)}\bigg).
\end{align*}
Then, for all values of $t$, $\rho$, $r$ in the given range, we have}
\begin{align*}
\sinh(r)\frac{\partial}{\partial\rho}F(t;\rho,r)+\sinh(\rho)\frac{\partial}{\partial r}F(t;\rho,r)
<0.
\end{align*}
\begin{proof}
We set
\begin{align*}
X:=\frac{\cosh(r/2)}{\cosh(\rho/2)},
\end{align*}
and compute
\begin{align*}
&\frac{\partial}{\partial\rho}F(t;\rho,r)=-\frac{re^{-r^{2}/(4t)}}{\sinh(r)}\frac{2\cosh^{2}(r/2)
\sinh(\rho/2)}{\cosh^{3}(\rho/2)}= \\
&-\frac{re^{-r^{2}/(4t)}}{\sinh(r)}2X^{2}\tanh(\rho/2)=-F(t;\rho,r)\frac{2X^{2}}{2X^{2}-1}
\tanh(\rho/2),
\end{align*}
and
\begin{align*}
\frac{\partial}{\partial r}F(t;\rho,r)=F(t;\rho,r)\bigg(\frac{1}{r}-\frac{r}{2t}+\frac{2X^{2}}
{2X^{2}-1}\tanh(r/2)-\frac{\cosh(r)}{\sinh(r)}\bigg).
\end{align*}
From this we deduce
\begin{align}
\notag
&\sinh(r)\frac{\partial}{\partial\rho}F(t;\rho,r)+\sinh(\rho)\frac{\partial}{\partial r}F(t;\rho,r)
= \\
\label{37}
&-F(t;\rho,r)\sinh(\rho)\bigg(\frac{r}{2t}+\frac{\cosh(r)}{\sinh(r)}-\frac{1}{r}\bigg)-F(t;\rho,r)
\frac{2X^{2}}{2X^{2}-1}h_{\rho}(r),
\end{align}
where
\begin{align*}
h_{\rho}(r):=\tanh(\rho/2)\sinh(r)-\sinh(\rho)\tanh(r/2).
\end{align*}
For $r>0$, we now have the estimate
\begin{align*}
\frac{r}{2t}+\frac{\cosh(r)}{\sinh(r)}-\frac{1}{r}>\frac{\cosh(r)}{\sinh(r)}-\frac{1}{r}=\frac{r
\cosh(r)-\sinh(r)}{r\sinh(r)}>0,
\end{align*}
using the power series expansions for $\cosh(r)$ and $\sinh(r)$. Next, we compute and
estimate for $r\geq\rho>0$
\begin{align*}
h_{\rho}'(r)&=\tanh(\rho/2)\cosh(r)-\frac{\sinh(\rho)}{2\cosh^{2}(r/2)} \\
&=\tanh(\rho/2)\big(2\cosh^{2}(r/2)-1\big)-\frac{\sinh(\rho/2)\cosh(\rho/2)}{\cosh^{2}
(r/2)} \\
&=\frac{\tanh(\rho/2)}{\cosh^{2}(r/2)}\big(2\cosh^{4}(r/2)-\cosh^{2}(r/2)-\cosh^{2}(\rho/2)
\big) \\
&\geq\frac{2\tanh(\rho/2)}{\cosh^{2}(r/2)}\big(\cosh^{4}(r/2)-\cosh^{2}(r/2)\big) \\
&=2\tanh(\rho/2)\big(\cosh^{2}(r/2)-1\big)>0. \\
\end{align*}
Since $h_{\rho}(\rho)=0$, this shows that $h_{\rho}(r)\geq 0$ for $r\geq\rho>0$. Recalling
\eqref{37} the claim of the lemma follows from the above estimates.
\end{proof}

\textbf{A.2. Proposition.} \emph{For any $t>0$, the heat kernel $K_{\mathbb{H}}^{(1)}(t;
\rho)$ for forms is strictly monotone decreasing for $\rho>0$.}
\begin{proof}
We will prove that $\partial/\partial_{\rho}K_{\mathbb{H}}^{(1)}(t;\rho)<0$ for $\rho>0$.
To simplify notations, we put
\begin{align*}
c(t):=\frac{\sqrt{2}e^{-t/4}}{(4\pi t)^{3/2}}.
\end{align*}
In the notation of Lemma A.1, we then have, using integration by parts,
\begin{align*}
K_{\mathbb{H}}^{(1)}(t;\rho)&=c(t)\int\limits_{\rho}^{\infty}F(t;\rho,r)\frac{\sinh(r)}{\sqrt
{\cosh(r)-\cosh(\rho)}}\,\mathrm{d}r \\
&=-2c(t)\int\limits_{\rho}^{\infty}\frac{\partial}{\partial r}F(t;\rho,r)\sqrt{\cosh(r)-\cosh
(\rho)}\,\mathrm{d}r\,.
\end{align*}
We now apply the Leibniz rule of differentiation to write
\begin{align*}
\frac{\partial}{\partial\rho}K_{\mathbb{H}}^{(1)}(t;\rho)&=-2c(t)\int\limits_{\rho}^{\infty}
\frac{\partial^{2}}{\partial r\,\partial\rho}F(t;\rho,r)\sqrt{\cosh(r)-\cosh(\rho)}\,\mathrm{d}
r \\
&\hspace*{4mm}+c(t)\int\limits_{\rho}^{\infty}\frac{\partial}{\partial r}F(t;\rho,r)\frac{\sinh
(\rho)}{\sqrt{\cosh(r)-\cosh(\rho)}}\,\mathrm{d}r\,.
\end{align*}
Using integration by parts on the first term once again, yields the identity
\begin{align*}
\frac{\partial}{\partial\rho}K_{\mathbb{H}}^{(1)}(t;\rho)=c(t)\int\limits_{\rho}^{\infty}\bigg(
\sinh(r)\frac{\partial}{\partial\rho}F(t;\rho,r)+\sinh(\rho)\frac{\partial}{\partial r}F(t;\rho,r)
\bigg)\frac{\mathrm{d}r}{\sqrt{\cosh(r)-\cosh(\rho)}}\,.
\end{align*}
From Lemma A.1, we conclude that  $\partial/\partial_{\rho}K_{\mathbb{H}}^{(1)}
(t;\rho)<0$ for $\rho>0$, which proves the claim.
\end{proof}

\newpage

\noindent
Jay Jorgenson \\
Department of Mathematics \\
City College of New York \\
Convent Avenue at 138th Street \\
New York, NY 10031
U.S.A. \\
e-mail: jjorgenson@mindspring.com

\vspace{5mm}

\noindent
J\"urg Kramer \\
Institut f\"ur Mathematik \\
Humboldt-Universit\"at zu Berlin \\
Unter den Linden 6 \\
D-10099 Berlin \\
Germany \\
e-mail: kramer@math.hu-berlin.de

\end{document}